\documentclass[12pt]{amsart}
\allowdisplaybreaks[2]
\usepackage{amssymb}
\usepackage{graphicx}
\numberwithin{equation}{section}
\theoremstyle{plain}
\newtheorem{thm}{Theorem}[section]
\newtheorem{theorem}[thm]{Theorem}

\newtheorem{corollary}[thm]{Corollary}

\theoremstyle{definition}

\newtheorem{definition}[thm]{Definition}

\begin{document}
\title[A Logical Calculus To Intuitively And Logically Denote Number Systems]
{A Logical Calculus To Intuitively And Logically Denote Number Systems} %
\author[Pith Xie]{Pith Xie}
\address{Department of Information Science and Communication \\ Nanjing
University of Information Science and Technology \\ Nanjing, 210044, China} %
\curraddr{P.O.Box 383 \\ Gulou Post Office \\ Gulou District %
\\ Nanjing, 210008, China} %
\email{pith.xie@gmail.com} %
%\keywords{j}
\begin{abstract}
Simple continued fractions, base-b expansions, Dedekind cuts and Cauchy sequences are common %
notations for number systems. In this note, first, it is proven that both simple continued %
fractions and base-b expansions fail to denote real numbers and thus lack logic; second, it %
is shown that Dedekind cuts and Cauchy sequences fail to join in algebraical operations and %
thus lack intuition; third, we construct a logical calculus and deduce numbers to intuitively %
and logically denote number systems. %
\end{abstract}
\maketitle
\setcounter{tocdepth}{5} \setcounter{page}{1}
%\tableofcontents
%\newpage
%
\section{Introduction}

Number system is a set together with one or more operations. Any notation for number system %
has to denote both set and operations. The common notations for number systems are simple %
continued fractions, base-b expansions, Dedekind cuts and Cauchy sequences. %

In \cite{Ref1} and \cite{Ref2}, simple continued fractions and base-b expansions denote each %
number in number systems as a set of symbols. So both them denote number systems %
intuitively and join well in algebraical operations. In \cite{Ref3} and \cite{Ref4}, Dedekind %
cuts and Cauchy sequences introduce infinite rational numbers to denote an irrational number. %
So both them denote number systems logically and join well in logical deduction. %

In this note, first, it is proven that both simple continued fractions and base-b expansions %
fail to denote real numbers and thus lack logic; second, it is shown that Dedekind cuts and %
Cauchy sequences fail to intuitively join in algebraical operations and thus lack intuition. %

However, mathematical logic has sufficiency of intuition and logic. In \cite{Ref9}, formal %
language introduces producer ``$ \rightarrow $" to formalize intuitive language. In %
\cite{Ref10}, propositional logic introduces connectives such as ``$ \neg $", ``$ \wedge $", %
``$ \vee $", ``$ \Rightarrow / \rightarrow $" and ``$ \Leftrightarrow / \leftrightarrow $" to %
formalize logical deduction. Therefore, it is feasible to combine producer and connectives to %
deduce intuitive and logical notations for number systems. %

The paper is organized as follows. In Section 2, we study the most common notation for number %
system ---\!--- decimals, and prove that they fail to denote real numbers. In Section 3, by %
comparing those common notations for number systems, we show that intuitive simple continued %
fractions and base-b expansions lack logic while logical Dedekind cuts and Cauchy sequences %
lack intuition. In Section 4, we construct a logical calculus and deduce numbers to %
intuitively and logically denote number systems. %

\section{Decimals And Real Number System}\label{Sec_DEC}

In this section, we show a conceptual error in the proof to \cite[THEOREM 134]{Ref1}, and %
then correct \cite[THEOREM 134]{Ref1}. %

\begin{definition}\label{Den_LIM}
A sequence $ \{ x_n \} $ in a metric space $ (X,d) $ is a convergent sequence if there exists %
a point $ x \in X $ such that, for every $ \epsilon > 0 $, there exists an integer $ N $ such %
that $ d(x,x_n) < \epsilon $ for every integer $ n \geq N $. The point $ x $ is called the %
limit of the sequence $ \{ x_n \} $ and we write %
\begin{eqnarray}
x_n \rightarrow x
\end{eqnarray}
or
\begin{eqnarray}
\lim\limits_{n \to \infty} x_n = x.
\end{eqnarray}
\end{definition}

\begin{theorem}[{\cite[THEOREM 134]{Ref1}}]\label{Them_SUM}
Any positive number $ \xi $ may be expressed as a decimal %
\begin{eqnarray}
A_{1} A_{2} \cdots A_{s+1}.\ a_1 a_2 a_3 \cdots,
\end{eqnarray}
where $ 0 \leq A_{1} < 10, 0 \leq A_{2} < 10, \cdots , 0 \leq a_{n} < 10 $, %
not all A and a are 0, and an infinity of the $ a_n $ are less than 9. If $ \xi \geq 1 $, %
then $ A_{1} \geq 0 $. There is a (1,1) correspondence between the numbers and the decimals, %
and %
\begin{eqnarray}\label{Eqn_XI}
\xi = A_{1} \cdot 10^s + \cdots + A_{s+1} + \frac{a_1}{10} +
\frac{a_2}{10^2} + \cdots .
\end{eqnarray}
\end{theorem}

\begin{proof}
Let $ [\xi] $ be the integral part of $ \xi $. Then we write %
\begin{eqnarray}\label{9.1.1}
\xi = [\xi] + x = X + x,
\end{eqnarray}
where $ X $ is an integer and $ 0 \leq x < 1 $, and consider $ X $ and $ x $ separately. %

If $ X > 0 $ and $ 10^{s} \leq x < 10^{s+1} $, %
and $ A_1 $ and $ X_1 $ are the quotient and remainder when $ X $ is divided by $ 10^{s} $, %
then $ X = A_1 \cdot 10^{s} + X_1 $, where $ 0 < A_1 = [10^{-s}X] < 10 $, %
$ 0 \leq X_1 < 10^s $. %

Similarly %
\begin{eqnarray*}
X_1 = & A_2 \cdot 10^{s-1} + X_2 & (0 \leq A_2 < 10, 0 \leq X_2 < 10^{s-1}), \\
X_2 = & A_3 \cdot 10^{s-2} + X_3 & (0 \leq A_3 < 10, 0 \leq X_3 < 10^{s-2}), \\
\cdots & \cdots & \cdots \\
X_{s-1} = & A_s \cdot 10 + X_s & (0 \leq A_s < 10, 0 \leq X_s < 10), \\
X_s = & A_{s+1} & (0 \leq A_{s+1} < 10).
\end{eqnarray*}
Thus $ X $ may be expressed uniquely in the form %
\begin{eqnarray}
X = A_1 \cdot 10^s + A_2 \cdot 10^{s-1} + \cdots + A_s \cdot 10 +
A_{s+1},
\end{eqnarray}
where every $ A $ is one of 0, 1, 2, $ \cdots $, 9, and $ A_1 $ is not 0. We abbreviate %
this expression to %
\begin{eqnarray}\label{9.1.3}
X = A_1 A_2 \cdots A_s A_{s+1},
\end{eqnarray}
the ordinary representation of $ X $ in decimal notation. %

Passing to $ x $, we write %
\begin{eqnarray*}
& X = f_1 & (0 \leq f_1 < 1).
\end{eqnarray*}
We suppose that $ a_1 = [10 f_1] $, so that %
\begin{eqnarray*}
\frac{a_1}{10} \leq f_1 < \frac{a_1 + 1}{10};
\end{eqnarray*}
$ a_1 $ is one of 0, 1, 2, $ \cdots $, 9, and %
\begin{eqnarray*}
a_1 = [10 f_1], & 10 f_1 = a_1 + f_2 & (0 \leq f_2 < 1).
\end{eqnarray*}

Similarly, we define $ a_2, a_3, \cdots $ by %
\begin{eqnarray*}
a_2 = [10 f_2], & 10 f_2 = a_2 + f_3 & (0 \leq f_3 < 1), \\
a_3 = [10 f_3], & 10 f_3 = a_3 + f_4 & (0 \leq f_4 < 1), \\
\cdots & \cdots & \cdots
\end{eqnarray*}
Every $ a_n $ is one of 0, 1, 2, $ \cdots $, 9. Thus %
\begin{eqnarray}\label{9.1.4}
x = x_n + g_{n+1},
\end{eqnarray}
where
\begin{eqnarray}
x_n = \frac{a_1}{10} + \frac{a_2}{10^2} + \cdots + \frac{a_n}{10^n}, \label{9.1.5} \\
0 \leq g_{n+1} = \frac{f_{n+1}}{10^n} < \frac{1}{10^n}.
\end{eqnarray}

We thus define a decimal $ . a_1 a_2 a_3 \cdots a_n \cdots $ associated with $ x $. We call %
$ a_1, a_2, \cdots $ the first, second, $ \cdots $ \emph{digits} of the decimal. %

Since $ a_n < 10 $, the series %
\begin{eqnarray}\label{9.1.7}
\sum\limits_{1}^{\infty} \frac{a_n}{10^n}
\end{eqnarray}
is convergent; and since $ g_{n+1} \rightarrow 0 $, its sum is $ x $. We may therefore write %
\begin{eqnarray}\label{9.1.8}
x = .\ a_1 a_2 a_3 \cdots,
\end{eqnarray}
the right-hand side being an abbreviation for the series (\ref{9.1.7}). %

If $ f_{n+1} = 0 $ for some $ n $, \emph{i.e.} if $ 10^n x $ is an integer, then %
\begin{eqnarray*}
a_{n+1} = a_{n+2} = \cdots = 0.
\end{eqnarray*}

In this case we say that the decimal \emph{terminates}. Thus %
\begin{eqnarray*}
\frac{17}{400} = .0425000 \cdots,
\end{eqnarray*}
and we write simply $ \frac{17}{400} = .0425 $. %

It is plain that the decimal for $ x $ will terminate if and only if $ x $ is a rational %
fraction whose denominator is of the form $ 2^\alpha 5^\beta $. %

Since $ \frac{a_{n+1}}{10^{n+1}} + \frac{a_{n+2}}{10^{n+2}} + \cdots = g_{n+1} < %
\frac{1}{10^n} $ and $ \frac{9}{10^{n+1}} + \frac{9}{10^{n+2}} + \cdots = %
\frac{9}{10^{n+1}(1-\frac{1}{10})} = \frac{1}{10^n} $, it is impossible that every $ a_n $ %
from a certain point on should be 9. With this reservation, every possible sequence $ (a_n) $ %
will arise from some $ x $. We define $ x $ as the sum of the series (\ref{9.1.7}), and %
$ x_n $ and $ g_{n+1} $ as in (\ref{9.1.4}) and (\ref{9.1.5}). Then $ g_{n+1} < 10^{-n} $ for %
every $ n $, and $ x $ yields the sequence required. %

Finally, if
\begin{eqnarray}\label{9.1.9}
\sum\limits_{1}^{\infty} \frac{a_n}{10^n} = \sum\limits_{1}^{\infty} \frac{b_n}{10^n}, %
\end{eqnarray}
and the $ b_n $ satisfy the conditions already imposed on the $ a_n $, then $ a_n = b_n $ for %
every $ n $. For if not, let $ a_N $ and $ b_N $ be the first pair which differ, so that %
$ |a_N - b_N| \geq 1 $. Then %
\begin{eqnarray*}
\left| \sum\limits_{1}^{\infty} \frac{a_n}{10^n} - \sum\limits_{1}^{\infty} \frac{b_n}{10^n} %
\right| \geq \frac{1}{10^N} - \sum\limits_{N+1}^{\infty} \frac{|a_n - b_n|}{10^n} \geq %
\frac{1}{10^N} - \sum\limits_{N+1}^{\infty} \frac{9}{10^n} = 0. %
\end{eqnarray*}

This contradicts (\ref{9.1.9}) unless there is equality. If there is equality, then all of %
$ a_{N+1} - b_{N+1}, a_{N+2} - b_{N+2}, \cdots $ must have the same sign and the absolute %
value 9. But then either $ a_n = 9 $ and $ b_n = 0 $ for $ n > N $, or else $ a_n = 0 $ and %
$ b_n = 9 $, and we have seen that each of these alternatives is impossible. Hence %
$ a_n = b_n $ for all $ n $. In other words, different decimals correspond to different %
numbers. %

We now combine (\ref{9.1.1}), (\ref{9.1.3}), and (\ref{9.1.8}) in the form %
\begin{eqnarray}
\xi = X + x = A_1 A_2 \cdots A_s A_{s+1}.\ a_1 a_2 a_3 \cdots; %
\end{eqnarray}
and the claim follows. %
\end{proof}

According to Definition \ref{Den_LIM}, the series (\ref{9.1.7}) converges to the limit $ x $. %
For an infinite sequence, however, its limit may not equal its $ \omega-th $ number for any %
infinite number $ \omega $. %

1. $ \omega $ is a \emph{transfinite cardinal number}\cite{Ref5}. Since the equalities and order on %
the fractions including transfinite cardinal numbers have not been defined, the equation %
$ g_{\omega+1} = \frac{f_{\omega+1}}{10^\omega} = 0 $ cannot be derived from %
given premises for any $ \omega $. %

2. $ \omega $ is an \emph{infinite superreal number}\cite{Ref6} or an %
\emph{infinite surreal number}\cite{Ref7}. Since the infinitesimal %
$ g_{\omega+1} = \frac{f_{\omega+1}}{10^\omega} > 0 $ holds for every $ \omega $, the %
equation $ g_{\omega+1} = \frac{f_{\omega+1}}{10^\omega} = 0 $ cannot be derived from given %
premises for any $ \omega $. %

In summary, the equation $ x = x_\omega + g_{\omega+1} $ cannot derives %
$ x = x_\omega $ for any infinite number $ \omega $. Thus, (\ref{9.1.8}) cannot be %
derived from given premises. %

In fact, the proof to \cite[THEOREM 134]{Ref1} confuses the limit and the $ \omega-th $ %
number of the same infinite sequence for some infinite number $ \omega $. %

According to the arguments above, we correct \cite[THEOREM 134]{Ref1} as follows. %
\begin{theorem}
Any positive number $ \xi $ may be expressed as a limit of an infinite decimal sequence %
\begin{eqnarray}
\lim\limits_{n \to \infty} A_{1} A_{2} \cdots A_{s+1}.\ a_1 a_2 a_3 \cdots a_n, %
\end{eqnarray}

where $ 0 \leq A_{1} < 10, 0 \leq A_{2} < 10, \cdots , 0 \leq a_{n} < 10 $, %
not all A and a are 0, and an infinity of the $ a_n $ are less than 9. If $ \xi \geq 1 $, %
then $ A_{1} \geq 0 $. There is a (1,1) correspondence between the numbers and the limits of %
infinite decimal sequences, and %
\begin{eqnarray}
\xi = A_{1} \cdot 10^s + \cdots + A_{s+1} + \lim\limits_{n \to \infty} \sum \frac{a_n}{10^n}. %
\end{eqnarray}
\end{theorem}

\section{Common Notations For Number Systems}

\subsection{Intuitive Notations}
Simple continued fractions and base-b expansions construct intuitive symbols to denote number %
systems. They join well in algebraical operations and thus have sufficiency of intuition. %

\begin{definition}
A finite continued fraction is a function %
\begin{eqnarray}\label{10.1.1}
a_0 + \cfrac{1}{a_1 + \cfrac{1}{a_2 + \cfrac{1}{\begin{matrix} a_3 + & \cdots \\
& + \cfrac{1}{a_N} \end{matrix}}}} %
\end{eqnarray}
of $ N+1 $ variables %
\begin{eqnarray}
a_0, a_1, \cdots, a_n, \cdots, a_N, %
\end{eqnarray}
which is called finite simple continued fraction when $ a_0, a_1, \cdots, a_N $ are integers %
such that $ a_n > 0 $ for all $ n \geq 1 $. %
\end{definition}

Finite simple continued fractions can be written in a compact abbreviated notation as %
\begin{eqnarray}
[a_0, a_1, a_2, \cdots, a_N]. %
\end{eqnarray}

\begin{definition}
If $ a_0, a_1, a_2, \cdots, a_n, \cdots $ is a sequence of integers such that $ a_n > 0 $ for %
all $ n \geq 1 $, then the notation %
\begin{eqnarray}
[a_0, a_1, a_2, \cdots] %
\end{eqnarray}
denotes an infinite simple continued fraction. %
\end{definition}

\begin{theorem}[{\cite[THEOREM 149]{Ref1}}]\label{The_CF1}
If $ p_{n} $ and $ q_{n} $ are defined by %
\begin{eqnarray}
p_{0}=a_{0}, & \ p_{1}=a_{1}a_{0}+1, & \ p_{n}=a_{n}p_{n-1}+p_{n-2} \ (2 \leq n \leq N), \\
q_{0}=1, & \ q_{1}=a_{1}, & \ q_{n}=a_{n}q_{n-1}+q_{n-2} \ (2 \leq n
\leq N),
\end{eqnarray}

then
\begin{eqnarray}
[ a_0, a_1, \ldots , a_n ] = \frac{p_{n}}{q_{n}}.
\end{eqnarray}
\end{theorem}

Theorem \ref{The_CF1} can be specialized for finite simple continued fractions as follows: %

\begin{theorem}\label{The_CF2}
$ \{ a_0, a_1, \cdots, a_n \} $ is an integer sequence. If $ p_{n} $ and $ q_{n} $ are %
defined by %
\begin{eqnarray}
p_{0}=a_{0}, & \ p_{1}=a_{1}a_{0}+1, & \ p_{n}=a_{n}p_{n-1}+p_{n-2} \ (2 \leq n \leq N), \\
q_{0}=1, & \ q_{1}=a_{1}, & \ q_{n}=a_{n}q_{n-1}+q_{n-2} \ (2 \leq n
\leq N),
\end{eqnarray}

then
\begin{eqnarray}
[ a_0, a_1, \ldots , a_n ] = \frac{p_{n}}{q_{n}}.
\end{eqnarray}
\end{theorem}

Theorem \ref{The_CF2} can directly derive such a corollary as follows: %
\begin{corollary}\label{Coy_CF2}
Any finite simple continued fraction can be represented by a rational number. %
\end{corollary}

\begin{theorem}[{\cite[THEOREM 161]{Ref1}}]
Any rational number can be represented by a finite simple continued fraction. %
\end{theorem}

According to Corollary \ref{Coy_CF2} and \cite[THEOREM 161]{Ref1}, finite simple continued %
fractions are equivalent to rational numbers. %

\begin{theorem}[{\cite[THEOREM 161]{Ref1}}]
Any rational number can be represented by a finite simple continued fraction. %
\end{theorem}

\begin{theorem}[{\cite[THEOREM 170]{Ref1}}]
Every irrational number can be expressed in just one way as an infinite simple continued %
fraction. %
\end{theorem}

\begin{proof}
We call %
\begin{eqnarray}
a'_n = [a_n, a_{n+1}, \cdots]
\end{eqnarray}
the n-th complete quotient of the continued fraction $ x = [a_0, a_1, \cdots] $. %

Clearly %
\begin{eqnarray*}
a'_n & = & \lim\limits_{N \to \infty} [a_n, a_{n+1}, \cdots, a_N] \\
& = & a_n + \lim\limits_{N \to \infty} \frac{1}{[a_{n+1}, \cdots, a_N]} \\
& = & a_n + \frac{1}{a'_{n+1}},
\end{eqnarray*}
and in particular $ x = a'_0 = a_0 + \frac{1}{a'_1} $. %

Also $ a'_n > a_n, a'_{n+1} > a_{n+1} > 0, 0 < \frac{1}{a'_{n+1}} < 1 $; and so %
$ a_n = [a'_n] $, the integral part of $ a'_n $. %

Let $ x $ be any real number, and let $ a_0 = [x] $. Then %
\begin{eqnarray*}
x = a_0 + \xi_0, \ \ \ 0 \leq \xi_0 < 1.
\end{eqnarray*}
If $ \xi_0 \neq 0 $, we can write %
\begin{eqnarray*}
\frac{1}{\xi_0} = a'_1, \ \ \ [a'_1] = a_1, \ \ \ a'_1 = a_1 + \xi_1, \ \ \ 0 \leq \xi_1 < 1. %
\end{eqnarray*}
If $ \xi_1 \neq 0 $, we can write %
\begin{eqnarray*}
\frac{1}{\xi_1} = a'_2 = a_2 + \xi_2, \ \ \ 0 \leq \xi_2 < 1, %
\end{eqnarray*}
and so on. Also $ a'_n = 1/\xi_{n-1} > 1 $, and so $ a_n \geq 1 $, for $ n \geq 1 $. Thus, %
\begin{eqnarray}
x = [a_0, a'_1] = \left[ a_0, a_1 + \frac{1}{a'_2} \right] = [a_0, a_1, a'_2] = %
[a_0, a_1, a_2, a'_3] = \cdots, %
\end{eqnarray}
where $ a_0, a_1, \cdots $ are integers and %
\begin{eqnarray}
a_1 > 0, \ \ \ a_2 > 0, \cdots. %
\end{eqnarray}

The system of equations %
\begin{eqnarray*}
x = & a_0 + \xi_0 & (0 \leq \xi_0 < 1), \\
\frac{1}{\xi_0} = & a'_1 = a_1 + \xi_1 & (0 \leq \xi_1 < 1), \\
\frac{1}{\xi_1} = & a'_2 = a_2 + \xi_2 & (0 \leq \xi_2 < 1), \\
\cdots & \cdots & \cdots
\end{eqnarray*}
is known as the \emph{continued fraction algorithm}. The algorithm continues so long as %
$ \xi_n \neq 0 $. If we eventually reach a value of $ n $, say $ N $, for which %
$ \xi_n = 0 $, the algorithm terminates and %
\begin{eqnarray}
x = [a_0, a_1, a_2, \cdots, a_N]. %
\end{eqnarray}
In this case $ x $ is represented by a simple continued fraction, and is rational. %

If $ x $ is an integer, then $ \xi_0 = 0 $ and $ x = a_0 $. If $ x $ is not integral, then %
\begin{eqnarray*}
x = \frac{h}{k},
\end{eqnarray*}
where $ h $ and $ k $ are integers and $ k > 1 $. Since %
\begin{eqnarray*}
\frac{h}{k} = a_0 + \xi_0, \ \ \ h = a_0 + \xi_0 k,
\end{eqnarray*}
$ a_0 $ is the quotient, and $ k_1 = \xi_0 k $ the remainder, when $ h $ is divided by $ k $. %

If $ \xi_0 \neq 0 $, then %
\begin{eqnarray}
a'_1 = \frac{1}{\xi_0} = \frac{k}{k_1}
\end{eqnarray}
and %
\begin{eqnarray*}
\frac{k}{k_1} = a_1 + \xi_1, \ \ \ k = a_1 k_1 + \xi_0 k_1;
\end{eqnarray*}
thus $ a_1 $ is the quotient, and $ k_2 = \xi_1 k_1 $ the remainder, when $ k $ is divided by %
$ k_1 $. We thus obtain a series of equations %
\begin{eqnarray*}
h = a_0 k + k_1, \ \ \ k = a_1 k_1 + k_2, \ \ \ k_1 = a_2 k_2 + k_3, \ \ \ \cdots %
\end{eqnarray*}
continuing so long as $ \xi_n \neq 0 $, or, what is the same thing, so long as %
$ k_{n+1} \neq 0 $. %

The non-negative integers $ k, k_1, k_2, \cdots $ form a strictly decreasing sequence, and so %
$ k_{N+1} = 0 $ for some $ N $. It follows that $ \xi_N = 0 $ for some $ N $, and that the %
continued fraction algorithm terminates. This proves \cite[THEOREM 161]{Ref1}. %

The system of equations %
\begin{eqnarray*}
h = & a_0 k + k_1 & (0 < k_1 < k), \\
k = & a_1 k_1 + k_2 & (0 < k_2 < k_1), \\
\cdots & \cdots & \cdots \\
k_{N-2} = & a_{N-1} k_{N-1} + k_{N} & (0 < k_{N} < k_{N-1}), \\
k_{N-1} = & a_{N} k_{N} &
\end{eqnarray*}
is known as \emph{Euclid's algorithm}.

If $ x $ is irrational the continued fraction algorithm cannot terminate. Hence it defines an %
infinite sequence of integers %
\begin{eqnarray}
a_0, a_1, a_2, \cdots,
\end{eqnarray}
and as before %
\begin{eqnarray}
x = [a_0, a'_1] = [a_0, a_1, a'_2] = \cdots = [a_0, a_1, a_2, \cdots, a_n, a'_{n+1}], %
\end{eqnarray}
where $ a'_{n+1} = a_{n+1} + \frac{1}{a'_{n+2}} > a_{n+1} $. %
Hence
\begin{eqnarray}
x = a'_0 = \frac{a'_1 a_0 +1}{a'_1} = \cdots = \frac{a'_{n+1}p_n + %
p_{n-1}}{a'_{n+1}q_n + q_{n-1}}, %
\end{eqnarray}
and so %
\begin{eqnarray}
x - \frac{p_n}{q_n} = \frac{p_{n-1}q_n - p_n q_{n-1}}{q_n(a'_{n+1}q_n + q_{n-1})} = %
\frac{(-1)^{n}}{q_n(a'_{n+1}q_n + q_{n-1})}, \\
|x - \frac{p_n}{q_n}| < \frac{1}{q_n(a_{n+1}q_n + q_{n-1})} = \frac{1}{q_n q_{n+1}} \leq %
\frac{1}{n(n+1)} \rightarrow 0, \label{For_CFI}
\end{eqnarray}
when $ n \rightarrow \infty $. Thus %
\begin{eqnarray}\label{For_CFE}
x = \lim\limits_{n \to \infty} \frac{p_n}{q_n} = [a_0, a_1, \cdots, a_n, \cdots], %
\end{eqnarray}
and the algorithm leads to the continued fraction whose value is $ x $. %
\end{proof}

In Section \ref{Sec_DEC}, we have proven that the limit of an infinite sequence may not equal %
the $ \omega-th $ number of the same infinite sequence for any infinite number $ \omega $. %

1. $ \omega $ is a transfinite cardinal number. Since the equalities and order on the %
fractions including transfinite cardinal numbers have not been defined, the inequality %
$ |x - \frac{p_\omega}{q_\omega}| < \frac{1}{\omega(\omega+1)} $ does not %
hold for any $ \omega $. %

2. $ \omega $ is an infinite superreal number or an infinite surreal number. Since the %
infinitesimal $ \frac{1}{\omega(\omega+1)} > 0 $ holds for every $ \omega $, %
the equation $ |x - \frac{p_\omega}{q_\omega}| = 0 $ or %
$ x = \frac{p_\omega}{q_\omega} $ cannot be derived from given premises for any %
$ \omega $. %

In summary, the inequality (\ref{For_CFI}) cannot derives (\ref{For_CFE}). In fact, %
(\ref{For_CFI}) only derives $ x = \lim\limits_{n \to \infty} \frac{p_n}{q_n} = %
\lim\limits_{n \to \infty} [a_0, a_1, \cdots, a_n] $. %

According to the arguments above, we correct \cite[THEOREM 170]{Ref1} as follows. %
\begin{theorem}
Every irrational number can be expressed in just one way as a limit of an infinite simple %
continued fraction sequence. %
\end{theorem}

According to \cite[\S BF.2]{Ref2}, we can define base-b expansions as follows: %

\begin{definition}\label{Den_BBE}
Base-b expansion is an expression of number as follows. %
\begin{eqnarray}
c_{n}b^{n} + c_{n-1}b^{n-1} \cdots + c_{2}b^{2} + c_{1}b^{1} +
c_{0}b^{0} + d_{1}b^{-1} + d_{2}b^{-2} \cdots + d_{n}b^{-n},
\end{eqnarray}
where $ b $ represents the base, and $ c_{i} $ and $ d_{i} $ are place-value coefficients. %
The expansion would ordinarily be written without the plus signs and the powers of the base %
as follows: %
\begin{eqnarray}
c_{n}c_{n- 1} \cdots c_{2}c_{1}c_{0}\ .\ d_{1}d_{2} \cdots d_{n},
\end{eqnarray}
where $ b^{i} $ is implied by the place-value property of the system. %
\end{definition}

According to \ref{Den_BBE}, finite decimals are just base-10 expansions. %

\begin{definition}\label{Den_BVE}
Base-variable expansions are base-b expansions for every finite integer b greater %
than 1. %
\end{definition}

\begin{theorem}\label{Them_BVE}
Every base-variable expansion is equal to a rational number. %
\end{theorem}

\begin{proof}
According to Definition \ref{Den_BVE}, every base-variable expansion $ x $ must also be a %
base-b expansion. Then %
\begin{eqnarray}
x = \pm a_{n}\cdots a_{2}a_{1}a_{0}\ .\ a_{-1}a_{-2} \cdots a_{-n} .
\end{eqnarray}

According to Definition \ref{Den_BBE}, it follows that %
\begin{eqnarray}\label{For_BVE}
&\pm a_{n} \cdots a_{2}a_{1}a_{0}\ .\ a_{-1}a_{-2} \cdots a_{-n} =
\pm \frac{{\sum\limits_{i=-n}^{n} a_{i}b^{i+n}}}{b^{n}}.&
\end{eqnarray}

Since both digit $ 0 \leq a_{i} < b $ and $ b $ are integers, %
$ \pm \frac{{\sum\limits_{i=-n}^{n} a_{i}b^{i+n}}}{b^{n}} $ must be a rational number. %
So the claim follows. %
\end{proof}

\begin{theorem}[The Fundamental Theorem of Arithmetic]\label{them_arc}
Every natural number is either prime or can be uniquely factored as
a product of primes in a unique way. %
\end{theorem}

\begin{theorem}\label{Them_BBE}
Every base-b expansion for a constant b may be unequal to a rational number. %
\end{theorem}

\begin{proof}
According to the equation (\ref{For_BVE}), every base-b expansion for a constant b may be %
expressed as follows: %
\begin{eqnarray}
x = \pm \frac{{\sum\limits_{i=-n}^{n} a_{i}b^{i+n}}}{b^{n}}.
\end{eqnarray}

Since there exists infinite primes, there must exist a prime $ q $ such that $ (q,b)=1 $. %
Since $ q>1 $ and $ b \neq 0 $, it follows from Theorem \ref{them_arc} that for every %
$ 0 \leq a_i < b $ there exists %
\begin{eqnarray}
q \cdot \sum\limits_{i=-n}^{n} a_{i}b^{i+n} \neq b^{n}.
\end{eqnarray}
Hence %
\begin{eqnarray}
\pm \frac{\sum\limits_{i=-n}^{n} a_{i}b^{i+n}}{b^{n}} & \neq & \frac{1}{q}, %
\end{eqnarray}
which holds for every $ 0 \leq a_i < b $. So the claim follows. %
\end{proof}

From Theorem \ref{Them_BVE} and Theorem \ref{Them_BBE}, we can conclude such a corollary as %
follows: %
\begin{corollary}
Base-variable expansions are included in rational numbers. %
\end{corollary}

According to the arguments above, no algorithms can determine the equalities between infinite %
simple continued fractions or infinite base-variable expansions and real numbers. As to the %
limits of infinite simple continued fraction sequences and those of infinite base-variable %
expansion sequences, they belong to logical notations and will be discussed in the next %
section. %

In summary, both simple continued fractions and base-variable expansions lack logic and fail %
to denote real numbers. %

\subsection{Logical Notations}

According to Definition \ref{Den_LIM}, limit is based on infinite sequence. So the limits of %
infinite simple continued fraction sequences and those of infinite base-variable expansion %
sequences are also defined on infinite sequence. %

In 1872, Dedekind and Cantor invented Dedekind cuts and Cauchy sequences respectively to %
denote number systems. However, both Dedekind cuts and Cauchy sequences are based on %
rational number system. In 1889, Peano published a study giving an axiomatic approach to %
the natural numbers\cite{Ref8}. Peano Axioms can also be extended to define rational numbers. %
Then both Dedekind cuts and Cauchy sequences join well in logical deduction and thus have %
sufficiency of logic. %

In nature, Dedekind cuts and Cauchy sequences introduce infinite rational numbers to denote %
an irrational number. In Dedekind cuts, an irrational cut $ (A, B) $ is defined on two %
infinite rational sets $ A $ and $ B $. In Cauchy sequences, an irrational number is defined %
on an equivalence class of some infinite rational sequence. %

Although it is feasible to logically define algebraical operations on infinite sets or %
infinite sequences, it is impossible to intuitively execute these infinite algebraical %
operations in a finite period. So the limits of infinite simple continued fraction sequences %
and those of infinite base-variable expansion sequences lack intuition and fail to join in %
algebraical operations. For the same reason, both Dedekind cuts and Cauchy sequences lack %
intuition and fail to join in algebraical operations. %

\section{Logical Calculus}
Simple continued fractions and base-variable expansions fail to denote real numbers, while %
the limits of infinite simple continued fraction sequences, the limits of infinite %
base-variable expansion sequences, logical Dedekind cuts and Cauchy sequences fail to join in %
algebraical operations. In mathematical logic, logical calculus is a formal system to %
abstract and analyze the induction and deduction apart from specific meanings. In this %
section, however, we construct a logical calculus by virtue of formal language and deduce %
numbers to intuitively and logically denote number systems. The logical calculus not only %
denotes real numbers, but also allows them to join in algebraical operations. %

The introduction of formal language aims to use computer fast execute real number operations. %
For clarity, we will explain the logical calculus with natural language. %

In \cite{Ref9}, the producer ``$ \rightarrow $" substitutes the right permutations for the %
left permutations to produce new permutations. In \cite{Ref10}, the connectives ``$ \neg $", %
``$ \wedge $", ``$ \vee $", ``$ \Rightarrow $" and ``$ \Leftrightarrow $" stand for ``not", %
``and", ``or", ``implies" and ``if and only if" respectively. Here, the producer %
``$ \rightarrow $" is considered as a predicate symbol and embedded into logical calculus. %

\begin{definition}
$ \{ \Phi, \Psi \} $ is a logical calculus such that: %
\begin{eqnarray}
\label{4.1} && \Phi \{ \\
\label{4.2} && V \{ \emptyset, a, b \cdots \}, \\
\label{4.3} && C \{ \emptyset, 1, + \cdots \}, \\
\label{4.4} && P \{ \emptyset, \in, \subseteq, \rightarrow, |, =, < \cdots \}, \\
\label{4.5} && V \circ C \{ \emptyset, a, b \cdots, 1, + \cdots, aa, ab \cdots, a1, %
a+ \cdots, ba, bb \cdots, b1, b+ \cdots, \\
\notag && aaa, aab \cdots, aa1, aa+ \cdots, baa, bab \cdots, ba1, ba+ \cdots \}, \\
\label{4.6} && C \circ C \{ \emptyset, 1, + \cdots, 11, 1+ \cdots, 111, 11+ \cdots \}, \\
\label{4.7} && V \circ C \circ P \{ \emptyset, a, b \cdots, 1, + \cdots, \in, \subseteq %
\cdots, aa, ab \cdots, a1, a+ \cdots, a\in, a\subseteq \cdots, \\
\notag && ba, bb \cdots, b1, b+ \cdots, b\in, b\subseteq \cdots, aaa, aab \cdots, %
aa1, aa+ \cdots, aa\in, aa\subseteq \cdots, \\
\notag && baa, bab \cdots, ba1, ba+ \cdots, ba\in, ba\subseteq \cdots \}, \\
\label{4.8} && (\hat{a} \in V) \Leftrightarrow ((\hat{a} \equiv \emptyset) \vee %
(\hat{a} \equiv a) \vee (\hat{a} \equiv b) \cdots), \\
\label{4.9} && (\hat{a} \in C) \Leftrightarrow ((\hat{a} \equiv \emptyset) \vee %
(\hat{a} \equiv 1) \vee (\hat{a} \equiv +) \cdots), \\
\label{4.10} && (\hat{a} \in (V \circ C)) \Leftrightarrow ((\hat{a} \equiv \emptyset) \vee %
(\hat{a} \equiv a) \vee (\hat{a} \equiv b) \cdots \vee (\hat{a} \equiv 1) \cdots \vee %
(\hat{a} \equiv aa) \\
\notag && \vee (\hat{a} \equiv ab) \cdots \vee (\hat{a} \equiv a1) \cdots \vee %
(\hat{a} \equiv aaa) \vee (\hat{a} \equiv aab) \cdots \vee (\hat{a} \equiv aa1) \cdots), \\
\label{4.11} && (\hat{a} \in (C \circ C)) \Leftrightarrow ((\hat{a} \equiv \emptyset) \vee %
(\hat{a} \equiv 1) \vee (\hat{a} \equiv +) \cdots \vee (\hat{a} \equiv 11) \vee %
(\hat{a} \equiv 1+) \cdots \\
\notag && \vee (\hat{a} \equiv 111) \vee (\hat{a} \equiv 11+) \cdots), \\
\label{4.12} && (\hat{a} \in (V \circ C \circ P)) \Leftrightarrow %
((\hat{a} \equiv \emptyset) \vee (\hat{a} \equiv a) \vee (\hat{a} \equiv b) \cdots \vee %
(\hat{a} \equiv \in) \cdots \vee (\hat{a} \equiv aa) \\
\notag && \vee (\hat{a} \equiv ab) \cdots \vee (\hat{a} \equiv a\in) \cdots \vee %
(\hat{a} \equiv aaa) \vee (\hat{a} \equiv aab) \cdots \vee (\hat{a} \equiv aa\in) \cdots), \\
\label{4.13} && (\bar{a} \in (V \circ C)) \wedge (\bar{b} \in (V \circ C)) \wedge %
(\bar{c} \in (V \circ C)) \wedge (\bar{d} \in (V \circ C)) \wedge (\bar{e} \in (V \circ C)) \\
\notag && \wedge (\bar{f} \in (V \circ C)) \wedge (\bar{g} \in (V \circ C)) \wedge %
(\bar{h} \in (V \circ C)) \wedge (\bar{i} \in (V \circ C)) \wedge (\bar{j} \in (V \circ C)) \\
\notag && \cdots \wedge (\bar{\bar{a}} \in (V \circ C \circ P)) \wedge %
(\bar{\bar{b}} \in (V \circ C \circ P)) \wedge %
(\bar{\bar{c}} \in (V \circ C \circ P)) \cdots, \\
\label{4.14} && ((\bar{a} \subseteq \{\bar{b},\bar{c}\}) \Leftrightarrow %
((\bar{a} \subseteq \bar{b}) \vee (\bar{a} \subseteq \bar{c}))) \wedge \\
\notag && ((\bar{a} \subseteq \{\bar{b},\bar{c},\bar{d}\}) \Leftrightarrow %
((\bar{a} \subseteq \bar{b}) \vee (\bar{a} \subseteq \bar{c}) \vee %
(\bar{a} \subseteq \bar{d}))) \wedge \\
\notag && ((\bar{a} \subseteq \{\bar{b},\bar{c},\bar{d},\bar{e}\}) \Leftrightarrow %
((\bar{a} \subseteq \bar{b}) \vee (\bar{a} \subseteq \bar{c}) \vee %
(\bar{a} \subseteq \bar{d}) \vee (\bar{a} \subseteq \bar{e}))) \wedge \\
\notag && ((\bar{a} \subseteq \{\bar{b},\bar{c},\bar{d},\bar{e},\bar{f}\}) \Leftrightarrow %
((\bar{a} \subseteq \bar{b}) \vee (\bar{a} \subseteq \bar{c}) \vee %
(\bar{a} \subseteq \bar{d}) \vee (\bar{a} \subseteq \bar{e}) \vee %
(\bar{a} \subseteq \bar{f}))) \wedge \\
\notag && ((\bar{a} \subseteq \{\bar{b},\bar{c},\bar{d},\bar{e},\bar{f},\bar{g}\}) %
\Leftrightarrow ((\bar{a} \subseteq \bar{b}) \vee (\bar{a} \subseteq \bar{c}) \vee %
(\bar{a} \subseteq \bar{d}) \vee (\bar{a} \subseteq \bar{e}) \vee (\bar{a} \subseteq \bar{f}) \\
\notag && \vee (\bar{a} \subseteq \bar{g}))) \wedge \\
\notag && ((\bar{a} \subseteq \{\bar{b},\bar{c},\bar{d},\bar{e},\bar{f},\bar{g},\bar{h}\}) %
\Leftrightarrow ((\bar{a} \subseteq \bar{b}) \vee (\bar{a} \subseteq \bar{c}) \vee %
(\bar{a} \subseteq \bar{d}) \vee (\bar{a} \subseteq \bar{e}) \vee (\bar{a} \subseteq \bar{f}) \\
\notag && \vee (\bar{a} \subseteq \bar{g}) \vee (\bar{a} \subseteq \bar{h}))) \wedge \\
\notag && ((\bar{a} \subseteq %
\{\bar{b},\bar{c},\bar{d},\bar{e},\bar{f},\bar{g},\bar{h},\bar{i}\}) %
\Leftrightarrow ((\bar{a} \subseteq \bar{b}) \vee (\bar{a} \subseteq \bar{c}) \vee %
(\bar{a} \subseteq \bar{d}) \vee (\bar{a} \subseteq \bar{e}) \vee (\bar{a} \subseteq \bar{f}) \\
\notag && \vee (\bar{a} \subseteq \bar{g}) \vee (\bar{a} \subseteq \bar{h}) %
\vee (\bar{a} \subseteq \bar{i}))) \wedge \\
\notag && ((\bar{a} \subseteq %
\{\bar{b},\bar{c},\bar{d},\bar{e},\bar{f},\bar{g},\bar{h},\bar{i},\bar{j}\}) %
\Leftrightarrow ((\bar{a} \subseteq \bar{b}) \vee (\bar{a} \subseteq \bar{c}) \vee %
(\bar{a} \subseteq \bar{d}) \vee (\bar{a} \subseteq \bar{e}) \vee (\bar{a} \subseteq \bar{f}) \\
\notag && \vee (\bar{a} \subseteq \bar{g}) \vee (\bar{a} \subseteq \bar{h}) %
\vee (\bar{a} \subseteq \bar{i}) \vee (\bar{a} \subseteq \bar{j}))) \cdots \\
\notag && \}, \\
\label{4.15} && \Psi \{ \\
\label{4.16} && (\bar{a} \subseteq \bar{b}) \Leftrightarrow %
(\bar{b} = \bar{c}\bar{a}\bar{d}), \\
\label{4.17} && (\bar{a} \rightarrow \bar{b}\bar{c}\bar{d}) \wedge %
(\bar{c} \rightarrow \bar{e}) \Rightarrow (\bar{a} \rightarrow \bar{b}\bar{e}\bar{d}), \\
\label{4.18} && (\bar{\bar{a}} \rightarrow \bar{\bar{b}} | \bar{\bar{c}}) \Rightarrow %
((\bar{\bar{a}} \rightarrow \bar{\bar{b}}) \wedge (\bar{\bar{a}} \rightarrow \bar{\bar{c}})), \\
\label{4.19} && (\bar{\bar{a}} | \bar{\bar{b}} \rightarrow \bar{\bar{c}}) \Rightarrow %
((\bar{\bar{a}} \rightarrow \bar{\bar{c}}) \wedge (\bar{\bar{b}} \rightarrow \bar{\bar{c}})), \\
\label{4.20} && (\bar{a} < \bar{b}) \Rightarrow \neg(\bar{b} < \bar{a}), \\
\label{4.21} && (\bar{a} < \bar{b}) \Rightarrow \neg(\bar{a} = \bar{b}), \\
\label{4.22} && (\bar{a} < \bar{b}) \wedge (\bar{b} < \bar{c}) \Rightarrow %
(\bar{a} < \bar{c}), \\
\label{4.23} && (\bar{a} < \bar{b}) \wedge (\bar{a} \in (C \circ C)) \wedge %
(\bar{b} \in (C \circ C)) \Rightarrow (\bar{a} \wedge \bar{b}), \\
\label{4.24} && (\bar{a} < \bar{b}\bar{c}\bar{d}) \wedge (\bar{c} = \bar{e}) \Rightarrow %
(\bar{a} < \bar{b}\bar{e}\bar{d}), \\
\label{4.25} && (\bar{a}\bar{b}\bar{c} < \bar{d}) \wedge (\bar{b} = \bar{e}) \Rightarrow %
(\bar{a}\bar{e}\bar{c} < \bar{d}), \\
\label{4.26} && (\bar{a} < \bar{b}\bar{c}\bar{d}) \wedge (\bar{c} \rightarrow \bar{e}) \wedge %
\neg(\bar{c} \subseteq \{\bar{a},\bar{b},\bar{d}\}) \Rightarrow %
(\bar{a} < \bar{b}\bar{e}\bar{d}), \\
\label{4.27} && (\bar{a} < \bar{b}\bar{c}\bar{d}\bar{c}\bar{e}) \wedge %
(\bar{c} \rightarrow \bar{f}) \wedge %
\neg(\bar{c} \subseteq \{\bar{a},\bar{b},\bar{d},\bar{e}\}) \Rightarrow %
(\bar{a} < \bar{b}\bar{f}\bar{d}\bar{f}\bar{e}), \\
\label{4.28} && (\bar{a}\bar{b}\bar{c} < \bar{d}) \wedge (\bar{b} \rightarrow \bar{e}) \wedge %
\neg(\bar{b} \subseteq \{\bar{a},\bar{c},\bar{d}\}) \Rightarrow %
(\bar{a}\bar{e}\bar{c} < \bar{d}), \\
\label{4.29} && (\bar{a}\bar{b}\bar{c} < \bar{d}\bar{b}\bar{e}) \wedge %
(\bar{b} \rightarrow \bar{f}) \wedge %
\neg(\bar{b} \subseteq \{\bar{a},\bar{c},\bar{d},\bar{e}\}) \Rightarrow %
(\bar{a}\bar{f}\bar{c} < \bar{d}\bar{f}\bar{e}), \\
\label{4.30} && (\bar{a}\bar{b}\bar{c} < \bar{d}\bar{b}\bar{e}\bar{b}\bar{f}) \wedge %
(\bar{b} \rightarrow \bar{g}) \wedge %
\neg(\bar{b} \subseteq \{\bar{a},\bar{c},\bar{d},\bar{e},\bar{f}\}) \Rightarrow %
(\bar{a}\bar{g}\bar{c} < \bar{d}\bar{g}\bar{e}\bar{g}\bar{f}), \\
\label{4.31} && (\bar{a}\bar{b}\bar{c}\bar{b}\bar{d} < \bar{e}) \wedge %
(\bar{b} \rightarrow \bar{f}) \wedge %
\neg(\bar{b} \subseteq \{\bar{a},\bar{c},\bar{d},\bar{e}\}) \Rightarrow %
(\bar{a}\bar{f}\bar{c}\bar{f}\bar{d} < \bar{e}), \\
\label{4.32} && (\bar{a}\bar{b}\bar{c}\bar{b}\bar{d} < \bar{e}\bar{b}\bar{f}) \wedge %
(\bar{b} \rightarrow \bar{g}) \wedge %
\neg(\bar{b} \subseteq \{\bar{a},\bar{c},\bar{d},\bar{e},\bar{f}\}) \Rightarrow %
(\bar{a}\bar{g}\bar{c}\bar{g}\bar{d} < \bar{e}\bar{g}\bar{f}), \\
\label{4.33} && (\bar{a}\bar{b}\bar{c}\bar{b}\bar{d} < \bar{e}\bar{b}\bar{f}\bar{b}\bar{g}) %
\wedge (\bar{b} \rightarrow \bar{h}) \wedge %
\neg(\bar{b} \subseteq \{\bar{a},\bar{c},\bar{d},\bar{e},\bar{f},\bar{g}\}) \Rightarrow %
(\bar{a}\bar{h}\bar{c}\bar{h}\bar{d} < \bar{e}\bar{h}\bar{f}\bar{h}\bar{g}), \\
\label{4.34} && \bar{a} = \bar{a}, \\
\label{4.35} && (\bar{a} = \bar{b}) \Rightarrow (\bar{b} = \bar{a}), \\
\label{4.36} && (\bar{a} = \bar{b}) \Rightarrow \neg (\bar{a} < \bar{b}), \\
\label{4.37} && (\bar{a} = \bar{b}\bar{c}\bar{d}) \wedge (\bar{c} = \bar{e}) \Rightarrow %
(\bar{a} = \bar{b}\bar{e}\bar{d}), \\
\label{4.38} && (\bar{a}\bar{b}\bar{c}) \wedge (\bar{b} = \bar{d}) \Rightarrow %
(\bar{a}\bar{b}\bar{c} = \bar{a}\bar{d}\bar{c}), \\
\label{4.39} && (\bar{a} = \bar{b}\bar{c}\bar{d}) \wedge (\bar{c} \rightarrow \bar{e}) \wedge %
\neg(\bar{c} \subseteq \{\bar{a},\bar{b},\bar{d}\}) \Rightarrow %
(\bar{a} = \bar{b}\bar{e}\bar{d}), \\
\label{4.40} && (\bar{a} = \bar{b}\bar{c}\bar{d}\bar{c}\bar{e}) \wedge %
(\bar{c} \rightarrow \bar{f}) \wedge %
\neg(\bar{c} \subseteq \{\bar{a},\bar{b},\bar{d},\bar{e}\}) \Rightarrow %
(\bar{a} = \bar{b}\bar{f}\bar{d}\bar{f}\bar{e}), \\
\label{4.41} && (\bar{a}\bar{b}\bar{c} = \bar{d}\bar{b}\bar{e}) \wedge %
(\bar{b} \rightarrow \bar{f}) \wedge %
\neg(\bar{b} \subseteq \{\bar{a},\bar{c},\bar{d},\bar{e}\}) \Rightarrow %
(\bar{a}\bar{f}\bar{c} = \bar{d}\bar{f}\bar{e}), \\
\label{4.42} && (\bar{a}\bar{b}\bar{c} = \bar{d}\bar{b}\bar{e}\bar{b}\bar{f}) \wedge %
(\bar{b} \rightarrow \bar{g}) \wedge %
\neg(\bar{b} \subseteq \{\bar{a},\bar{c},\bar{d},\bar{e},\bar{f}\}) \Rightarrow %
(\bar{a}\bar{g}\bar{c} = \bar{d}\bar{g}\bar{e}\bar{g}\bar{f}), \\
\label{4.43} && (\bar{a}\bar{b}\bar{c}\bar{b}\bar{d} = \bar{e}\bar{b}\bar{f}\bar{b}\bar{g}) %
\wedge (\bar{b} \rightarrow \bar{h}) \wedge %
\neg(\bar{b} \subseteq \{\bar{a},\bar{c},\bar{d},\bar{e},\bar{f},\bar{g}\}) \Rightarrow %
(\bar{a}\bar{h}\bar{c}\bar{h}\bar{d} = \bar{e}\bar{h}\bar{f}\bar{h}\bar{g}) \\
\notag && \}. %
\end{eqnarray}
\end{definition}

First, we will explain the primitive symbols of the logical calculus $ \{ \Phi, \Psi \} $ %
with natural language. %

The symbols ``$ \{ $", ``$ \} $", ``$ , $", ``$ ( $", ``$ ) $" are punctuation. The symbol %
``$ \emptyset $" indicates emptiness. The symbol ``$ \cdots $" indicates an omission. %

(\ref{4.1}) denotes $ \Phi $ as a set of notations and particular axioms between $ \{ $ and %
$ \} $. Different logical calculus correspond to different notations and particular axioms. %

(\ref{4.2}) denotes $ V $ as a set of variables between $ \{ $ and $ \} $. %

(\ref{4.3}) denotes $ C $ as a set of constants between $ \{ $ and $ \} $. %

(\ref{4.4}) denotes $ P $ as a set of predicate symbols between $ \{ $ and $ \} $. %

(\ref{4.5}) denotes $ V \circ C $ as a set of concatenations between $ V $ and $ C $. %

(\ref{4.6}) denotes $ C \circ C $ as a set of concatenations between $ C $ and $ C $. %

(\ref{4.7}) denotes $ V \circ C \circ P $ as a set of concatenations among $ V $, $ C $ and %
$ P $. %

(\ref{4.8}) $ \sim $ (\ref{4.12}) define a set of axioms on the binary predicate symbol %
$ \in $. %

(\ref{4.13}) defines an axiom on new variables ranging over $ V \circ C $. %

(\ref{4.14}) defines an axiom on the binary predicate symbol $ \subseteq $. %

(\ref{4.15}) denotes $ \Psi $ as a set of general axioms between $ \{ $ and $ \} $. %
Different logical calculus correspond to the same general axioms. %

(\ref{4.16}) defines an axiom on the binary predicate symbol $ \subseteq $. %

(\ref{4.17}) defines an axiom on the binary predicate symbol $ \rightarrow $. %

(\ref{4.18}) $ \sim $ (\ref{4.19}) define a set of axioms on the binary predicate symbol $ | $. %

(\ref{4.20}) $ \sim $ (\ref{4.33}) define a set of axioms on the binary predicate symbol $ < $. %

(\ref{4.34}) $ \sim $ (\ref{4.43}) define a set of axioms on the binary predicate symbol $ = $. %

%%%%%%%%%%%%%%%%%%%%%%%%%%%%%%%%%%%%%%%%%%%%%%%%%%%%%%%%%%%%%

Then, we will prove that the logical calculus $ \{ \Phi, \Psi \} $ can deduce common number %
systems. %

\begin{definition}
In a logical calculus $ \{ \Phi, \Psi \} $, if $ \bar{a} \equiv true $, then $ \bar{a} $ is a %
number. %
\end{definition}

\begin{theorem}
\begin{eqnarray}
\notag && \textit{If} \ \Phi \{ \\
\label{4.44} && V \{ \emptyset, a, b \}, \\
\label{4.45} && C \{ \emptyset, 1, + \}, \\
\label{4.46} && P \{ \emptyset, \in, \subseteq, \rightarrow, |, =, < \}, \\
\label{4.47} && V \circ C \{ \emptyset, a, b \cdots, 1, + \cdots, aa, ab \cdots, a1, %
a+ \cdots, ba, bb \cdots, b1, b+ \cdots, \\
\notag && aaa, aab \cdots, aa1, aa+ \cdots, baa, bab \cdots, ba1, ba+ \cdots \}, \\
\label{4.48} && C \circ C \{ \emptyset, 1, + \cdots, 11, 1+ \cdots, 111, 11+ \cdots \}, \\
\label{4.49} && V \circ C \circ P \{ \emptyset, a, b \cdots, 1, + \cdots, \in, \subseteq %
\cdots, aa, ab \cdots, a1, a+ \cdots, a\in, a\subseteq \cdots, \\
\notag && ba, bb \cdots, b1, b+ \cdots, b\in, b\subseteq \cdots, aaa, aab \cdots, %
aa1, aa+ \cdots, aa\in, aa\subseteq \cdots, \\
\notag && baa, bab \cdots, ba1, ba+ \cdots, ba\in, ba\subseteq \cdots \}, \\
\label{4.50} && (\hat{a} \in V) \Leftrightarrow ((\hat{a} \equiv \emptyset) \vee %
(\hat{a} \equiv a) \vee (\hat{a} \equiv b) \cdots), \\
\label{4.51} && (\hat{a} \in C) \Leftrightarrow ((\hat{a} \equiv \emptyset) \vee %
(\hat{a} \equiv 1) \vee (\hat{a} \equiv +) \cdots), \\
\label{4.52} && (\hat{a} \in (V \circ C)) \Leftrightarrow ((\hat{a} \equiv \emptyset) \vee %
(\hat{a} \equiv a) \vee (\hat{a} \equiv b) \cdots \vee (\hat{a} \equiv 1) \cdots \vee %
(\hat{a} \equiv aa) \\
\notag && \vee (\hat{a} \equiv ab) \cdots \vee (\hat{a} \equiv a1) \cdots \vee %
(\hat{a} \equiv aaa) \vee (\hat{a} \equiv aab) \cdots \vee (\hat{a} \equiv aa1) \cdots), \\
\label{4.53} && (\hat{a} \in (C \circ C)) \Leftrightarrow ((\hat{a} \equiv \emptyset) \vee %
(\hat{a} \equiv 1) \vee (\hat{a} \equiv +) \cdots \vee (\hat{a} \equiv 11) \vee %
(\hat{a} \equiv 1+) \cdots \\
\notag && \vee (\hat{a} \equiv 111) \vee (\hat{a} \equiv 11+) \cdots), \\
\label{4.54} && (\hat{a} \in (V \circ C \circ P)) \Leftrightarrow %
((\hat{a} \equiv \emptyset) \vee (\hat{a} \equiv a) \vee (\hat{a} \equiv b) \cdots \vee %
(\hat{a} \equiv \in) \cdots \vee (\hat{a} \equiv aa) \\
\notag && \vee (\hat{a} \equiv ab) \cdots \vee (\hat{a} \equiv a\in) \cdots \vee %
(\hat{a} \equiv aaa) \vee (\hat{a} \equiv aab) \cdots \vee (\hat{a} \equiv aa\in) \cdots), \\
\label{4.55} && (\bar{a} \in (V \circ C)) \wedge (\bar{b} \in (V \circ C)) \wedge %
(\bar{c} \in (V \circ C)) \wedge (\bar{d} \in (V \circ C)) \wedge (\bar{e} \in (V \circ C)) \\
\notag && \wedge (\bar{f} \in (V \circ C)) \wedge (\bar{g} \in (V \circ C)) \wedge %
(\bar{h} \in (V \circ C)) \wedge (\bar{i} \in (V \circ C)) \wedge (\bar{j} \in (V \circ C)) \\
\notag && \wedge (\bar{\bar{a}} \in (V \circ C \circ P)) \wedge %
(\bar{\bar{b}} \in (V \circ C \circ P)) \wedge (\bar{\bar{c}} \in (V \circ C \circ P)), \\
\label{4.56} && ((\bar{a} \subseteq \{\bar{b},\bar{c}\}) \Leftrightarrow %
((\bar{a} \subseteq \bar{b}) \vee (\bar{a} \subseteq \bar{c}))) \wedge \\
\notag && ((\bar{a} \subseteq \{\bar{b},\bar{c},\bar{d}\}) \Leftrightarrow %
((\bar{a} \subseteq \bar{b}) \vee (\bar{a} \subseteq \bar{c}) \vee %
(\bar{a} \subseteq \bar{d}))) \wedge \\
\notag && ((\bar{a} \subseteq \{\bar{b},\bar{c},\bar{d},\bar{e}\}) \Leftrightarrow %
((\bar{a} \subseteq \bar{b}) \vee (\bar{a} \subseteq \bar{c}) \vee %
(\bar{a} \subseteq \bar{d}) \vee (\bar{a} \subseteq \bar{e}))) \wedge \\
\notag && ((\bar{a} \subseteq \{\bar{b},\bar{c},\bar{d},\bar{e},\bar{f}\}) \Leftrightarrow %
((\bar{a} \subseteq \bar{b}) \vee (\bar{a} \subseteq \bar{c}) \vee %
(\bar{a} \subseteq \bar{d}) \vee (\bar{a} \subseteq \bar{e}) \vee %
(\bar{a} \subseteq \bar{f}))) \wedge \\
\notag && ((\bar{a} \subseteq \{\bar{b},\bar{c},\bar{d},\bar{e},\bar{f},\bar{g}\}) %
\Leftrightarrow ((\bar{a} \subseteq \bar{b}) \vee (\bar{a} \subseteq \bar{c}) \vee %
(\bar{a} \subseteq \bar{d}) \vee (\bar{a} \subseteq \bar{e}) \vee (\bar{a} \subseteq \bar{f}) \\
\notag && \vee (\bar{a} \subseteq \bar{g}))) \wedge \\
\notag && ((\bar{a} \subseteq \{\bar{b},\bar{c},\bar{d},\bar{e},\bar{f},\bar{g},\bar{h}\}) %
\Leftrightarrow ((\bar{a} \subseteq \bar{b}) \vee (\bar{a} \subseteq \bar{c}) \vee %
(\bar{a} \subseteq \bar{d}) \vee (\bar{a} \subseteq \bar{e}) \vee (\bar{a} \subseteq \bar{f}) \\
\notag && \vee (\bar{a} \subseteq \bar{g}) \vee (\bar{a} \subseteq \bar{h}))) \wedge \\
\notag && ((\bar{a} \subseteq %
\{\bar{b},\bar{c},\bar{d},\bar{e},\bar{f},\bar{g},\bar{h},\bar{i}\}) %
\Leftrightarrow ((\bar{a} \subseteq \bar{b}) \vee (\bar{a} \subseteq \bar{c}) \vee %
(\bar{a} \subseteq \bar{d}) \vee (\bar{a} \subseteq \bar{e}) \vee (\bar{a} \subseteq \bar{f}) \\
\notag && \vee (\bar{a} \subseteq \bar{g}) \vee (\bar{a} \subseteq \bar{h}) %
\vee (\bar{a} \subseteq \bar{i}))) \wedge \\
\notag && ((\bar{a} \subseteq %
\{\bar{b},\bar{c},\bar{d},\bar{e},\bar{f},\bar{g},\bar{h},\bar{i},\bar{j}\}) %
\Leftrightarrow ((\bar{a} \subseteq \bar{b}) \vee (\bar{a} \subseteq \bar{c}) \vee %
(\bar{a} \subseteq \bar{d}) \vee (\bar{a} \subseteq \bar{e}) \vee (\bar{a} \subseteq \bar{f}) \\
\notag && \vee (\bar{a} \subseteq \bar{g}) \vee (\bar{a} \subseteq \bar{h}) %
\vee (\bar{a} \subseteq \bar{i}) \vee (\bar{a} \subseteq \bar{j}))), \\
\label{4.57} && a \rightarrow 1|1+a, \\
\label{4.58} && a < 1+a, \\
\label{4.59} && \bar{a} \wedge \bar{b} \Rightarrow (\bar{a}+\bar{b} = \bar{b}+\bar{a}) \\
\notag && \}, %
\end{eqnarray}

then $ N \{ \Phi, \Psi \} $ denotes natural number system. %
\end{theorem}

\begin{proof}
\begin{eqnarray*}
(A1) & (a \rightarrow 1|1+a) \Rightarrow (a \rightarrow 1) & by (\ref{4.57}),(\ref{4.18}) \\
(A2) & \Rightarrow (a \rightarrow 1+a) & by (\ref{4.19}) \\
(A3) & (a < 1+a) \Rightarrow (1 < 1+1) & by (\ref{4.58}),(A1),(\ref{4.29}) \\
(A4) & \Rightarrow 1 & by (\ref{4.23}) \\
(A5) & \Rightarrow (1+1) & by (A3),(\ref{4.23}) \\
(A6) & \Rightarrow (1+a < 1+1+a) & by (\ref{4.58}),(A2),(\ref{4.29}) \\
(A7) & \Rightarrow (1+1 < 1+1+1) & by (A1),(\ref{4.29}) \\
(A8) & \Rightarrow (1+1) & by (\ref{4.23}) \\
(A9) & \Rightarrow (1+1+1) & by (A7),(\ref{4.23}) \\
\vdots & \vdots & \vdots
\end{eqnarray*}

Then we deduce the numbers from $ N \{ \Phi, \Psi \} $: %
\begin{eqnarray*}
& 1, 1+1, 1+1+1, 1+1+1+1 \cdots &
\end{eqnarray*}
\begin{eqnarray*}
(B1) & 1+1 = 1+1 & by (A4),(\ref{4.59}) \\
(B2) & 1+1+1 = 1+1+1 & by (A4),(A5),(\ref{4.59}) \\
(B3) & 1+1+1+1 = 1+1+1+1 & by (A5),(\ref{4.59}) \\
\vdots & \vdots & \vdots
\end{eqnarray*}

Then we deduce the equalities on deducible numbers from $ N \{ \Phi, \Psi \} $: %
\begin{eqnarray*}
& 1+1 = 1+1, 1+1+1 = 1+1+1, 1+1+1+1 = 1+1+1+1 \cdots &
\end{eqnarray*}

The deducible numbers correspond to the natural numbers as follows: %
\begin{eqnarray*}
1 & \equiv & 1, \\
1+1 & \equiv & 2, \\
1+1+1 & \equiv & 3, \\
\vdots & \vdots & \vdots .
\end{eqnarray*}

The equalities on deducible numbers correspond to the addition in natural number system. %
So the claim follows. %
\end{proof}

%%%%%%%%%%%%%%%%%%%%%%%%%%%%%%%%%%%%%%%%%%%%%%%%%%%%%%%%%%%%%

\begin{theorem}
\begin{eqnarray}
\notag && \textit{If} \ \Phi \{ \\
\label{4.60} && V \{ \emptyset, a, b, c \}, \\
\label{4.61} && C \{ \emptyset, 1, +, [, ], - \}, \\
\label{4.62} && P \{ \emptyset, \in, \subseteq, \rightarrow, |, =, < \}, \\
\label{4.63} && V \circ C \{ \emptyset, a, b \cdots, 1, + \cdots, aa, ab \cdots, a1, %
a+ \cdots, ba, bb \cdots, b1, b+ \cdots, \\
\notag && aaa, aab \cdots, aa1, aa+ \cdots, baa, bab \cdots, ba1, ba+ \cdots \}, \\
\label{4.64} && C \circ C \{ \emptyset, 1, + \cdots, 11, 1+ \cdots, 111, 11+ \cdots \}, \\
\label{4.65} && V \circ C \circ P \{ \emptyset, a, b \cdots, 1, + \cdots, \in, \subseteq %
\cdots, aa, ab \cdots, a1, a+ \cdots, a\in, a\subseteq \cdots, \\
\notag && ba, bb \cdots, b1, b+ \cdots, b\in, b\subseteq \cdots, aaa, aab \cdots, %
aa1, aa+ \cdots, aa\in, aa\subseteq \cdots, \\
\notag && baa, bab \cdots, ba1, ba+ \cdots, ba\in, ba\subseteq \cdots \}, \\
\label{4.66} && (\hat{a} \in V) \Leftrightarrow ((\hat{a} \equiv \emptyset) \vee %
(\hat{a} \equiv a) \vee (\hat{a} \equiv b) \cdots), \\
\label{4.67} && (\hat{a} \in C) \Leftrightarrow ((\hat{a} \equiv \emptyset) \vee %
(\hat{a} \equiv 1) \vee (\hat{a} \equiv +) \cdots), \\
\label{4.68} && (\hat{a} \in (V \circ C)) \Leftrightarrow ((\hat{a} \equiv \emptyset) \vee %
(\hat{a} \equiv a) \vee (\hat{a} \equiv b) \cdots \vee (\hat{a} \equiv 1) \cdots \vee %
(\hat{a} \equiv aa) \\
\notag && \vee (\hat{a} \equiv ab) \cdots \vee (\hat{a} \equiv a1) \cdots \vee %
(\hat{a} \equiv aaa) \vee (\hat{a} \equiv aab) \cdots \vee (\hat{a} \equiv aa1) \cdots), \\
\label{4.69} && (\hat{a} \in (C \circ C)) \Leftrightarrow ((\hat{a} \equiv \emptyset) \vee %
(\hat{a} \equiv 1) \vee (\hat{a} \equiv +) \cdots \vee (\hat{a} \equiv 11) \vee %
(\hat{a} \equiv 1+) \cdots \\
\notag && \vee (\hat{a} \equiv 111) \vee (\hat{a} \equiv 11+) \cdots), \\
\label{4.70} && (\hat{a} \in (V \circ C \circ P)) \Leftrightarrow %
((\hat{a} \equiv \emptyset) \vee (\hat{a} \equiv a) \vee (\hat{a} \equiv b) \cdots \vee %
(\hat{a} \equiv \in) \cdots \vee (\hat{a} \equiv aa) \\
\notag && \vee (\hat{a} \equiv ab) \cdots \vee (\hat{a} \equiv a\in) \cdots \vee %
(\hat{a} \equiv aaa) \vee (\hat{a} \equiv aab) \cdots \vee (\hat{a} \equiv aa\in) \cdots), \\
\label{4.71} && (\bar{a} \in (V \circ C)) \wedge (\bar{b} \in (V \circ C)) \wedge %
(\bar{c} \in (V \circ C)) \wedge (\bar{d} \in (V \circ C)) \wedge (\bar{e} \in (V \circ C)) \\
\notag && \wedge (\bar{f} \in (V \circ C)) \wedge (\bar{g} \in (V \circ C)) \wedge %
(\bar{h} \in (V \circ C)) \wedge (\bar{i} \in (V \circ C)) \wedge (\bar{j} \in (V \circ C)) \\
\notag && \wedge (\bar{\bar{a}} \in (V \circ C \circ P)) \wedge %
(\bar{\bar{b}} \in (V \circ C \circ P)) \wedge (\bar{\bar{c}} \in (V \circ C \circ P)), \\
\label{4.72} && ((\bar{a} \subseteq \{\bar{b},\bar{c}\}) \Leftrightarrow %
((\bar{a} \subseteq \bar{b}) \vee (\bar{a} \subseteq \bar{c}))) \wedge \\
\notag && ((\bar{a} \subseteq \{\bar{b},\bar{c},\bar{d}\}) \Leftrightarrow %
((\bar{a} \subseteq \bar{b}) \vee (\bar{a} \subseteq \bar{c}) \vee %
(\bar{a} \subseteq \bar{d}))) \wedge \\
\notag && ((\bar{a} \subseteq \{\bar{b},\bar{c},\bar{d},\bar{e}\}) \Leftrightarrow %
((\bar{a} \subseteq \bar{b}) \vee (\bar{a} \subseteq \bar{c}) \vee %
(\bar{a} \subseteq \bar{d}) \vee (\bar{a} \subseteq \bar{e}))) \wedge \\
\notag && ((\bar{a} \subseteq \{\bar{b},\bar{c},\bar{d},\bar{e},\bar{f}\}) \Leftrightarrow %
((\bar{a} \subseteq \bar{b}) \vee (\bar{a} \subseteq \bar{c}) \vee %
(\bar{a} \subseteq \bar{d}) \vee (\bar{a} \subseteq \bar{e}) \vee %
(\bar{a} \subseteq \bar{f}))) \wedge \\
\notag && ((\bar{a} \subseteq \{\bar{b},\bar{c},\bar{d},\bar{e},\bar{f},\bar{g}\}) %
\Leftrightarrow ((\bar{a} \subseteq \bar{b}) \vee (\bar{a} \subseteq \bar{c}) \vee %
(\bar{a} \subseteq \bar{d}) \vee (\bar{a} \subseteq \bar{e}) \vee (\bar{a} \subseteq \bar{f}) \\
\notag && \vee (\bar{a} \subseteq \bar{g}))) \wedge \\
\notag && ((\bar{a} \subseteq \{\bar{b},\bar{c},\bar{d},\bar{e},\bar{f},\bar{g},\bar{h}\}) %
\Leftrightarrow ((\bar{a} \subseteq \bar{b}) \vee (\bar{a} \subseteq \bar{c}) \vee %
(\bar{a} \subseteq \bar{d}) \vee (\bar{a} \subseteq \bar{e}) \vee (\bar{a} \subseteq \bar{f}) \\
\notag && \vee (\bar{a} \subseteq \bar{g}) \vee (\bar{a} \subseteq \bar{h}))) \wedge \\
\notag && ((\bar{a} \subseteq %
\{\bar{b},\bar{c},\bar{d},\bar{e},\bar{f},\bar{g},\bar{h},\bar{i}\}) %
\Leftrightarrow ((\bar{a} \subseteq \bar{b}) \vee (\bar{a} \subseteq \bar{c}) \vee %
(\bar{a} \subseteq \bar{d}) \vee (\bar{a} \subseteq \bar{e}) \vee (\bar{a} \subseteq \bar{f}) \\
\notag && \vee (\bar{a} \subseteq \bar{g}) \vee (\bar{a} \subseteq \bar{h}) %
\vee (\bar{a} \subseteq \bar{i}))) \wedge \\
\notag && ((\bar{a} \subseteq %
\{\bar{b},\bar{c},\bar{d},\bar{e},\bar{f},\bar{g},\bar{h},\bar{i},\bar{j}\}) %
\Leftrightarrow ((\bar{a} \subseteq \bar{b}) \vee (\bar{a} \subseteq \bar{c}) \vee %
(\bar{a} \subseteq \bar{d}) \vee (\bar{a} \subseteq \bar{e}) \vee (\bar{a} \subseteq \bar{f}) \\
\notag && \vee (\bar{a} \subseteq \bar{g}) \vee (\bar{a} \subseteq \bar{h}) %
\vee (\bar{a} \subseteq \bar{i}) \vee (\bar{a} \subseteq \bar{j}))), \\
\label{4.73} && a \rightarrow 1|[aba], \\
\label{4.74} && b|c \rightarrow +|-, \\
\label{4.75} && a < [1+a], \\
\label{4.76} && \bar{a} \wedge \bar{b} \wedge \bar{c} \Rightarrow %
([\bar{a}b\bar{b}c\bar{c}] = [[\bar{a}b\bar{b}]c\bar{c}]), \\
\label{4.77} && \bar{a} \Rightarrow ([\bar{a}-\bar{a}] = [1-1]), \\
\label{4.78} && \bar{a} \wedge \bar{b} \Rightarrow ([\bar{a}+\bar{b}] = [\bar{b}+\bar{a}]), \\
\label{4.79} && \bar{a} \wedge \bar{b} \Rightarrow ([[\bar{a}-\bar{b}]+\bar{b}] = \bar{a}), \\
\label{4.80} && \bar{a} \wedge \bar{b} \wedge \bar{c} \Rightarrow %
([\bar{a}-\bar{b}+\bar{c}] = [\bar{a}+\bar{c}-\bar{b}]), \\
\label{4.81} && \bar{a} \wedge \bar{b} \wedge \bar{c} \Rightarrow %
([\bar{a}+[\bar{b}+\bar{c}]] = [[\bar{a}+\bar{b}]+\bar{c}]), \\
\label{4.82} && \bar{a} \wedge \bar{b} \wedge \bar{c} \Rightarrow %
([\bar{a}+[\bar{b}-\bar{c}]] = [[\bar{a}+\bar{b}]-\bar{c}]), \\
\label{4.83} && \bar{a} \wedge \bar{b} \wedge \bar{c} \Rightarrow %
([\bar{a}-[\bar{b}+\bar{c}]] = [[\bar{a}-\bar{b}]-\bar{c}]), \\
\label{4.84} && \bar{a} \wedge \bar{b} \wedge \bar{c} \Rightarrow %
([\bar{a}-[\bar{b}-\bar{c}]] = [[\bar{a}-\bar{b}]+\bar{c}] \\
\notag && \}, %
\end{eqnarray}

then $ Z \{ \Phi, \Psi \} $ denotes integral number system. %
\end{theorem}

\begin{proof}
\begin{eqnarray*}
(A1) & (a \rightarrow 1|[aba]) \Rightarrow (a \rightarrow 1) & by (\ref{4.73}),(\ref{4.18}) \\
(A2) & \Rightarrow (a \rightarrow [aba]) & by (\ref{4.73}),(\ref{4.18}) \\
(A3) & \Rightarrow (a \rightarrow [1b1]) & by (A2),(A1),(\ref{4.17}) \\
(A4) & (b|c \rightarrow +|-) \Rightarrow (b|c \rightarrow -) & by (\ref{4.74}),(\ref{4.18}) \\
(A5) & \Rightarrow (b \rightarrow -) & by (\ref{4.19}) \\
(A6) & \Rightarrow (a \rightarrow [1-1]) & by (A3),(A5),(\ref{4.17}) \\
(A7) & (a < [1+a]) \Rightarrow (1 < [1+1]) & by (\ref{4.75}),(A1),(\ref{4.29}) \\
(A8) & \Rightarrow 1 & by (\ref{4.23}) \\
(A9) & \Rightarrow [1+1] & \ by (A7),(\ref{4.23}) \\
(A10) & (a < [1+a]) \Rightarrow ([1-1] < [1+[1-1]]) & by (\ref{4.75}),(A6),(\ref{4.29}) \\
(A11) & \Rightarrow ([1-1] < [[1-1]+1]) & by (\ref{4.78}),(\ref{4.24}) \\
(A12) & \Rightarrow ([1-1] < 1) & by (\ref{4.79}),(\ref{4.24}) \\
(A13) & \Rightarrow [1-1] & by (\ref{4.23}) \\
\vdots & \vdots & \vdots
\end{eqnarray*}

Then we deduce the numbers from $ Z \{ \Phi, \Psi \} $: %
\begin{eqnarray*}
& [1-1], 1, [1-1-1], [1+1], [1+1+1] \cdots & %
\end{eqnarray*}
\begin{eqnarray*}
(B1) & [1+1] = [1+1] & by (A8),(\ref{4.78}) \\
(B2) & [1+[1+1]] = [[1+1]+1] & by (A8),(A9),(\ref{4.78}) \\
(B3) & [1+[1-1]] = [[1-1]+1] & by (A8),(A13),(\ref{4.78}) \\
(B4) & [[1+1]+[1-1]] = [[1-1]+[1+1]] & by (A9),(A13),(\ref{4.78}) \\
\vdots & \vdots & \vdots
\end{eqnarray*}

Then we deduce the equalities on deducible numbers from $ Z \{ \Phi, \Psi \} $: %
\begin{eqnarray*}
& [1+[1+1]] = [[1+1]+1], [1+[1-1]] = [[1-1]+1] \cdots & %
\end{eqnarray*}

The deducible numbers correspond to the integral numbers as follows: %
\begin{eqnarray*}
\vdots & \vdots & \vdots, \\
\ [1-1-[1+1]] & \equiv & -2, \\
\ [1-1-1] & \equiv & -1, \\
\ [1-1] & \equiv & 0, \\
\ 1 & \equiv & 1, \\
\ [1+1] & \equiv & 2, \\
\ [1+1+1] & \equiv & 3, \\
\ [1+1+1+1] & \equiv & 4, \\
\vdots & \vdots & \vdots .
\end{eqnarray*}

The equalities on deducible numbers correspond to the addition and subtraction in integral %
number system. So the claim follows. %
\end{proof}

%%%%%%%%%%%%%%%%%%%%%%%%%%%%%%%%%%%%%%%%%%%%%%%%%%%%%%%%%%%%%

\begin{theorem}
\begin{eqnarray}
\notag && \textit{If} \ \Phi \{ \\
\label{4.85} && V \{ \emptyset, a, b, c, d \}, \\
\label{4.86} && C \{ \emptyset, 1, +, [, ], - \}, \\
\label{4.87} && P \{ \emptyset, \in, \subseteq, \rightarrow, |, =, < \}, \\
\label{4.88} && V \circ C \{ \emptyset, a, b \cdots, 1, + \cdots, aa, ab \cdots, a1, %
a+ \cdots, ba, bb \cdots, b1, b+ \cdots, \\
\notag && aaa, aab \cdots, aa1, aa+ \cdots, baa, bab \cdots, ba1, ba+ \cdots \}, \\
\label{4.89} && C \circ C \{ \emptyset, 1, + \cdots, 11, 1+ \cdots, 111, 11+ \cdots \}, \\
\label{4.90} && V \circ C \circ P \{ \emptyset, a, b \cdots, 1, + \cdots, \in, \subseteq %
\cdots, aa, ab \cdots, a1, a+ \cdots, a\in, a\subseteq \cdots, \\
\notag && ba, bb \cdots, b1, b+ \cdots, b\in, b\subseteq \cdots, aaa, aab \cdots, %
aa1, aa+ \cdots, aa\in, aa\subseteq \cdots, \\
\notag && baa, bab \cdots, ba1, ba+ \cdots, ba\in, ba\subseteq \cdots \}, \\
\label{4.91} && (\hat{a} \in V) \Leftrightarrow ((\hat{a} \equiv \emptyset) \vee %
(\hat{a} \equiv a) \vee (\hat{a} \equiv b) \cdots), \\
\label{4.92} && (\hat{a} \in C) \Leftrightarrow ((\hat{a} \equiv \emptyset) \vee %
(\hat{a} \equiv 1) \vee (\hat{a} \equiv +) \cdots), \\
\label{4.93} && (\hat{a} \in (V \circ C)) \Leftrightarrow ((\hat{a} \equiv \emptyset) \vee %
(\hat{a} \equiv a) \vee (\hat{a} \equiv b) \cdots \vee (\hat{a} \equiv 1) \cdots \vee %
(\hat{a} \equiv aa) \\
\notag && \vee (\hat{a} \equiv ab) \cdots \vee (\hat{a} \equiv a1) \cdots \vee %
(\hat{a} \equiv aaa) \vee (\hat{a} \equiv aab) \cdots \vee (\hat{a} \equiv aa1) \cdots), \\
\label{4.94} && (\hat{a} \in (C \circ C)) \Leftrightarrow ((\hat{a} \equiv \emptyset) \vee %
(\hat{a} \equiv 1) \vee (\hat{a} \equiv +) \cdots \vee (\hat{a} \equiv 11) \vee %
(\hat{a} \equiv 1+) \cdots \\
\notag && \vee (\hat{a} \equiv 111) \vee (\hat{a} \equiv 11+) \cdots), \\
\label{4.95} && (\hat{a} \in (V \circ C \circ P)) \Leftrightarrow %
((\hat{a} \equiv \emptyset) \vee (\hat{a} \equiv a) \vee (\hat{a} \equiv b) \cdots \vee %
(\hat{a} \equiv \in) \cdots \vee (\hat{a} \equiv aa) \\
\notag && \vee (\hat{a} \equiv ab) \cdots \vee (\hat{a} \equiv a\in) \cdots \vee %
(\hat{a} \equiv aaa) \vee (\hat{a} \equiv aab) \cdots \vee (\hat{a} \equiv aa\in) \cdots), \\
\label{4.96} && (\bar{a} \in (V \circ C)) \wedge (\bar{b} \in (V \circ C)) \wedge %
(\bar{c} \in (V \circ C)) \wedge (\bar{d} \in (V \circ C)) \wedge (\bar{e} \in (V \circ C)) \\
\notag && \wedge (\bar{f} \in (V \circ C)) \wedge (\bar{g} \in (V \circ C)) \wedge %
(\bar{h} \in (V \circ C)) \wedge (\bar{i} \in (V \circ C)) \wedge \\
\notag && (\bar{j} \in (V \circ C)) \wedge (\bar{\bar{a}} \in (V \circ C \circ P)) \wedge %
(\bar{\bar{b}} \in (V \circ C \circ P)) \wedge (\bar{\bar{c}} \in (V \circ C \circ P)), \\
\label{4.97} && ((\bar{a} \subseteq \{\bar{b},\bar{c}\}) \Leftrightarrow %
((\bar{a} \subseteq \bar{b}) \vee (\bar{a} \subseteq \bar{c}))) \wedge \\
\notag && ((\bar{a} \subseteq \{\bar{b},\bar{c},\bar{d}\}) \Leftrightarrow %
((\bar{a} \subseteq \bar{b}) \vee (\bar{a} \subseteq \bar{c}) \vee %
(\bar{a} \subseteq \bar{d}))) \wedge \\
\notag && ((\bar{a} \subseteq \{\bar{b},\bar{c},\bar{d},\bar{e}\}) \Leftrightarrow %
((\bar{a} \subseteq \bar{b}) \vee (\bar{a} \subseteq \bar{c}) \vee %
(\bar{a} \subseteq \bar{d}) \vee (\bar{a} \subseteq \bar{e}))) \wedge \\
\notag && ((\bar{a} \subseteq \{\bar{b},\bar{c},\bar{d},\bar{e},\bar{f}\}) \Leftrightarrow %
((\bar{a} \subseteq \bar{b}) \vee (\bar{a} \subseteq \bar{c}) \vee %
(\bar{a} \subseteq \bar{d}) \vee (\bar{a} \subseteq \bar{e}) \vee %
(\bar{a} \subseteq \bar{f}))) \wedge \\
\notag && ((\bar{a} \subseteq \{\bar{b},\bar{c},\bar{d},\bar{e},\bar{f},\bar{g}\}) %
\Leftrightarrow ((\bar{a} \subseteq \bar{b}) \vee (\bar{a} \subseteq \bar{c}) \vee %
(\bar{a} \subseteq \bar{d}) \vee (\bar{a} \subseteq \bar{e}) \vee (\bar{a} \subseteq \bar{f}) \\
\notag && \vee (\bar{a} \subseteq \bar{g}))) \wedge \\
\notag && ((\bar{a} \subseteq \{\bar{b},\bar{c},\bar{d},\bar{e},\bar{f},\bar{g},\bar{h}\}) %
\Leftrightarrow ((\bar{a} \subseteq \bar{b}) \vee (\bar{a} \subseteq \bar{c}) \vee %
(\bar{a} \subseteq \bar{d}) \vee (\bar{a} \subseteq \bar{e}) \vee (\bar{a} \subseteq \bar{f}) \\
\notag && \vee (\bar{a} \subseteq \bar{g}) \vee (\bar{a} \subseteq \bar{h}))) \wedge \\
\notag && ((\bar{a} \subseteq %
\{\bar{b},\bar{c},\bar{d},\bar{e},\bar{f},\bar{g},\bar{h},\bar{i}\}) %
\Leftrightarrow ((\bar{a} \subseteq \bar{b}) \vee (\bar{a} \subseteq \bar{c}) \vee %
(\bar{a} \subseteq \bar{d}) \vee (\bar{a} \subseteq \bar{e}) \vee (\bar{a} \subseteq \bar{f}) \\
\notag && \vee (\bar{a} \subseteq \bar{g}) \vee (\bar{a} \subseteq \bar{h}) %
\vee (\bar{a} \subseteq \bar{i}))) \wedge \\
\notag && ((\bar{a} \subseteq %
\{\bar{b},\bar{c},\bar{d},\bar{e},\bar{f},\bar{g},\bar{h},\bar{i},\bar{j}\}) %
\Leftrightarrow ((\bar{a} \subseteq \bar{b}) \vee (\bar{a} \subseteq \bar{c}) \vee %
(\bar{a} \subseteq \bar{d}) \vee (\bar{a} \subseteq \bar{e}) \vee \\
\notag && (\bar{a} \subseteq \bar{f}) \vee (\bar{a} \subseteq \bar{g}) \vee %
(\bar{a} \subseteq \bar{h}) \vee (\bar{a} \subseteq \bar{i}) \vee %
(\bar{a} \subseteq \bar{j}))), \\
\label{4.98} && a \rightarrow 1|[aba], \\
\label{4.99} && b \rightarrow +|-, \\
\label{4.100} && c|d \rightarrow b|++|--, \\
\label{4.101} && a < [1+a], \\
\label{4.102} && (\bar{a} < \bar{b}) \wedge \bar{c} \Rightarrow %
([\bar{a}+\bar{c}] < [\bar{b}+\bar{c}]), \\
\label{4.103} && (\bar{a} < \bar{b}) \wedge \bar{c} \Rightarrow %
([\bar{c}-\bar{b}] < [\bar{c}-\bar{a}]), \\
\label{4.104} && ([1-1] < \bar{a}) \wedge (\bar{a} < \bar{b}) \wedge ([1-1] < \bar{c}) %
\Rightarrow ([\bar{a}--\bar{c}] < [\bar{b}--\bar{c}]), \\
\label{4.105} && \bar{a} \wedge \bar{b} \wedge \bar{c} \Rightarrow %
([\bar{a}c\bar{b}d\bar{c}] = [[\bar{a}c\bar{b}]d\bar{c}]), \\
\label{4.106} && \bar{a} \Rightarrow ([\bar{a}-\bar{a}] = [1-1]), \\
\label{4.107} && \bar{a} \wedge \bar{b} \Rightarrow ([\bar{a}+\bar{b}] = [\bar{b}+\bar{a}]), \\
\label{4.108} && \bar{a} \wedge \bar{b} \Rightarrow ([[\bar{a}-\bar{b}]+\bar{b}] = \bar{a}), \\
\label{4.109} && \bar{a} \wedge \bar{b} \wedge \bar{c} \Rightarrow %
([\bar{a}-\bar{b}+\bar{c}] = [\bar{a}+\bar{c}-\bar{b}]), \\
\label{4.110} && \bar{a} \wedge \bar{b} \wedge \bar{c} \Rightarrow %
([\bar{a}+[\bar{b}+\bar{c}]] = [[\bar{a}+\bar{b}]+\bar{c}]), \\
\label{4.111} && \bar{a} \wedge \bar{b} \wedge \bar{c} \Rightarrow %
([\bar{a}+[\bar{b}-\bar{c}]] = [[\bar{a}+\bar{b}]-\bar{c}]), \\
\label{4.112} && \bar{a} \wedge \bar{b} \wedge \bar{c} \Rightarrow %
([\bar{a}-[\bar{b}+\bar{c}]] = [[\bar{a}-\bar{b}]-\bar{c}]), \\
\label{4.113} && \bar{a} \wedge \bar{b} \wedge \bar{c} \Rightarrow %
([\bar{a}-[\bar{b}-\bar{c}]] = [[\bar{a}-\bar{b}]+\bar{c}]), \\
\label{4.114} && \bar{a} \Rightarrow ([\bar{a}++1] = \bar{a}), \\
\label{4.115} && \neg (\bar{a} = [1-1]) \Rightarrow ([\bar{a}--\bar{a}] = 1), \\
\label{4.116} && \bar{a} \wedge \bar{b} \Rightarrow %
([\bar{a}++\bar{b}] = [\bar{b}++\bar{a}]), \\
\label{4.117} && \bar{a} \wedge \bar{b} \wedge \bar{c} \Rightarrow %
([\bar{a}++[\bar{b}+\bar{c}]] = [[\bar{a}++\bar{b}]+[\bar{a}++\bar{c}]]), \\
\label{4.118} && \bar{a} \wedge \bar{b} \wedge \bar{c} \Rightarrow %
([\bar{a}++[\bar{b}-\bar{c}]] = [[\bar{a}++\bar{b}]-[\bar{a}++\bar{c}]]), \\
\label{4.119} && \bar{a} \wedge \bar{b} \wedge \bar{c} \Rightarrow %
([\bar{a}++[\bar{b}++\bar{c}]] = [[\bar{a}++\bar{b}]++\bar{c}]), \\
\label{4.120} && \bar{a} \wedge \neg (\bar{b} = [1-1]) \Rightarrow %
([\bar{a}--\bar{b}++\bar{b}] = \bar{a}), \\
\label{4.121} && \bar{a} \wedge \bar{b} \wedge \neg (\bar{c} = [1-1]) \Rightarrow %
(([\bar{a}--\bar{c}++\bar{b}] = [\bar{a}++\bar{b}--\bar{c}]) \wedge %
([[\bar{a}+\bar{b}]--\bar{c}] = \\
\notag && [[\bar{a}--\bar{c}]+[\bar{b}--\bar{c}]]) \wedge %
([[\bar{a}-\bar{b}]--\bar{c}] = [[\bar{a}--\bar{c}]-[\bar{b}--\bar{c}]]) \wedge \\
\notag && ([\bar{a}++[\bar{b}--\bar{c}]] = [[\bar{a}++\bar{b}]--\bar{c}])), \\
\label{4.122} && \bar{a} \wedge \neg ([\bar{b} = [1-1]) \wedge \neg (\bar{c} = [1-1]]) %
\Rightarrow (([\bar{a}--[\bar{b}++\bar{c}]] = [[\bar{a}--\bar{b}]--\bar{c}]) \wedge \\
\notag && ([\bar{a}--[\bar{b}--\bar{c}]] = [[\bar{a}--\bar{b}]++\bar{c}])) \\
\notag && \}, %
\end{eqnarray}

then $ Q \{ \Phi, \Psi \} $ denotes rational number system. %
\end{theorem}

\begin{proof}
\begin{eqnarray*}
(A1) & (a \rightarrow 1|[aba]) \Rightarrow (a \rightarrow 1) & by (\ref{4.98}),(\ref{4.18}) \\
(A2) & \Rightarrow (a \rightarrow [aba]) & by (\ref{4.98}),(\ref{4.18}) \\
(A3) & \Rightarrow (a \rightarrow [[aba]ba]) & by (A2),(\ref{4.17}) \\
(A4) & \Rightarrow (a \rightarrow [[1-1]-1]) & by (A1),(\ref{4.99}) \\
(A5) & (a < [1+a]) \Rightarrow (1 < [1+1]) & by (\ref{4.101}),(A1),(\ref{4.29}) \\
(A6) & \Rightarrow 1 & by (\ref{4.23}) \\
(A7) & \Rightarrow [1+1] & \ by (A5),(\ref{4.23}) \\
(A8) & \Rightarrow ([aba] < [1+[aba]]) & by (A5),(A2),(\ref{4.29}) \\
(A9) & \Rightarrow ([1+1] < [1+[1+1]]) & by (A1),(\ref{4.99}),(\ref{4.29}) \\
(A10) & \Rightarrow ([1+1] < [[1+1]+1]) & by (\ref{4.107}),(\ref{4.24}) \\
(A11) & \Rightarrow ([1+1] < [1+1+1]) & by (\ref{4.105}) \\
(A12) & \Rightarrow ([1-1] < [1+[1-1]]) & by (A8),(A1),(\ref{4.99}),(\ref{4.29}) \\
(A13) & \Rightarrow ([1-1] < [[1-1]+1]) & by (\ref{4.107}),(\ref{4.24}) \\
(A14) & \Rightarrow ([1-1] < 1) & by (\ref{4.108}) \\
(A15) & \Rightarrow ([1-1] < [1+1+1]) & by (A14),(A5),(A11),(\ref{4.22}) \\
(A16) & \Rightarrow ([1--[1+1+1]] < [[1+1]--[1+1+1]]) & by (A14),(A5),(A15),(\ref{4.104}) \\
(A17) & \Rightarrow [1--[1+1+1]] & by (\ref{4.23}) \\
(A18) & \Rightarrow [[1+1]--[1+1+1]] & by (A16),(\ref{4.23}) \\
\vdots & \vdots & \vdots
\end{eqnarray*}

Then we deduce the numbers from $ Q \{ \Phi, \Psi \} $: %
\begin{eqnarray*}
& [1-1], [[1-1]-[1--[1+1]]], [1--[1+1]], 1 \cdots & %
\end{eqnarray*}
\begin{eqnarray*}
(B1) & [1+[1--[1+1+1]]] = [[1--[1+1+1]]+1] & by (A6),(A17),(\ref{4.107}) \\
(B2) & [1++[1--[1+1+1]]] = [[1--[1+1+1]]++1] & by (A6),(A17),(\ref{4.116}) \\
\vdots & \vdots & \vdots
\end{eqnarray*}

Then we deduce the equalities on deducible numbers from $ Q \{ \Phi, \Psi \} $: %
\begin{eqnarray*}
& [1+1] = [1+1], [1++[1--[1+1+1]]] = [[1--[1+1+1]]++1] \cdots & %
\end{eqnarray*}

The deducible numbers correspond to the rational numbers as follows: %
\begin{eqnarray*}
\vdots & \vdots & \vdots, \\
\ [1-1-1] & \equiv & -1, \\
\vdots & \vdots & \vdots, \\
\ [1-1-[1--[1+1]]] & \equiv & -\frac{1}{2}, \\
\vdots & \vdots & \vdots, \\
\ [1-1] & \equiv & 0, \\
\vdots & \vdots & \vdots, \\
\ [1-[1--[1+1]]] & \equiv & \frac{1}{2}, \\
\vdots & \vdots & \vdots, \\
\ 1 & \equiv & 1, \\
\vdots & \vdots & \vdots, \\
\ [1+[1--[1+1]]] & \equiv & \frac{3}{2}, \\
\vdots & \vdots & \vdots, \\
\ [1+1] & \equiv & 2, \\
\vdots & \vdots & \vdots, \\
\ [1+1+[1--[1+1]]] & \equiv & \frac{5}{2}, \\
\vdots & \vdots & \vdots, \\
\ [1+1+1] & \equiv & 3, \\
\vdots & \vdots & \vdots .
\end{eqnarray*}

The equalities on deducible numbers correspond to the addition, subtraction, multiplication, %
division in rational number system. So the claim follows. %
\end{proof}

%%%%%%%%%%%%%%%%%%%%%%%%%%%%%%%%%%%%%%%%%%%%%%%%%%%%%%%%%%%%%

\begin{definition}
Real number system is a logical calculus $ R \{ \Phi, \Psi \} $ such that: %
\begin{eqnarray}
\notag && \Phi \{ \\
\label{4.123} && V \{ \emptyset, a, b, c, d, e, f, g, h, i, j, k, l \}, \\
\label{4.124} && C \{ \emptyset, 1, +, [, ], -, /, \top, \bot, \underline{\ \ } \}, \\
\label{4.125} && P \{ \emptyset, \in, \subseteq, \rightarrow, |, =, <, \| \}, \\
\label{4.126} && V \circ C \{ \emptyset, a, b \cdots, 1, + \cdots, aa, ab \cdots, a1, %
a+ \cdots, ba, bb \cdots, b1, b+ \cdots, \\
\notag && aaa, aab \cdots, aa1, aa+ \cdots, baa, bab \cdots, ba1, ba+ \cdots \}, \\
\label{4.127} && C \circ C \{ \emptyset, 1, + \cdots, 11, 1+ \cdots, 111, 11+ \cdots \}, \\
\label{4.128} && V \circ C \circ P \{ \emptyset, a, b \cdots, 1, + \cdots, \in, \subseteq %
\cdots, aa, ab \cdots, a1, a+ \cdots, a\in, a\subseteq \cdots, \\
\notag && ba, bb \cdots, b1, b+ \cdots, b\in, b\subseteq \cdots, aaa, aab \cdots, %
aa1, aa+ \cdots, aa\in, aa\subseteq \cdots, \\
\notag && baa, bab \cdots, ba1, ba+ \cdots, ba\in, ba\subseteq \cdots \}, \\
\label{4.129} && (\hat{a} \in V) \Leftrightarrow ((\hat{a} \equiv \emptyset) \vee %
(\hat{a} \equiv a) \vee (\hat{a} \equiv b) \cdots), \\
\label{4.130} && (\hat{a} \in C) \Leftrightarrow ((\hat{a} \equiv \emptyset) \vee %
(\hat{a} \equiv 1) \vee (\hat{a} \equiv +) \cdots), \\
\label{4.131} && (\hat{a} \in (V \circ C)) \Leftrightarrow ((\hat{a} \equiv \emptyset) \vee %
(\hat{a} \equiv a) \vee (\hat{a} \equiv b) \cdots \vee (\hat{a} \equiv 1) \cdots \vee %
(\hat{a} \equiv aa) \\
\notag && \vee (\hat{a} \equiv ab) \cdots \vee (\hat{a} \equiv a1) \cdots \vee %
(\hat{a} \equiv aaa) \vee (\hat{a} \equiv aab) \cdots \vee (\hat{a} \equiv aa1) \cdots), \\
\label{4.132} && (\hat{a} \in (C \circ C)) \Leftrightarrow ((\hat{a} \equiv \emptyset) \vee %
(\hat{a} \equiv 1) \vee (\hat{a} \equiv +) \cdots \vee (\hat{a} \equiv 11) \vee %
(\hat{a} \equiv 1+) \cdots \\
\notag && \vee (\hat{a} \equiv 111) \vee (\hat{a} \equiv 11+) \cdots), \\
\label{4.133} && (\hat{a} \in (V \circ C \circ P)) \Leftrightarrow %
((\hat{a} \equiv \emptyset) \vee (\hat{a} \equiv a) \vee (\hat{a} \equiv b) \cdots \vee %
(\hat{a} \equiv \in) \cdots \vee (\hat{a} \equiv aa) \\
\notag && \vee (\hat{a} \equiv ab) \cdots \vee (\hat{a} \equiv a\in) \cdots \vee %
(\hat{a} \equiv aaa) \vee (\hat{a} \equiv aab) \cdots \vee (\hat{a} \equiv aa\in) \cdots), \\
\label{4.134} && (\bar{a} \in (V \circ C)) \wedge (\bar{b} \in (V \circ C)) \wedge %
(\bar{c} \in (V \circ C)) \wedge (\bar{d} \in (V \circ C)) \wedge (\bar{e} \in (V \circ C)) \\
\notag && \wedge (\bar{f} \in (V \circ C)) \wedge (\bar{g} \in (V \circ C)) \wedge %
(\bar{h} \in (V \circ C)) \wedge (\bar{i} \in (V \circ C)) \wedge \\
\notag && (\bar{j} \in (V \circ C)) \wedge (\bar{k} \in (V \circ C)) \wedge %
(\bar{l} \in (V \circ C)) \wedge (\bar{\bar{a}} \in (V \circ C \circ P)) \wedge \\
\notag && (\bar{\bar{b}} \in (V \circ C \circ P)) \wedge %
(\bar{\bar{c}} \in (V \circ C \circ P)), \\
\label{4.135} && ((\bar{a} \subseteq \{\bar{b},\bar{c}\}) \Leftrightarrow %
((\bar{a} \subseteq \bar{b}) \vee (\bar{a} \subseteq \bar{c}))) \wedge \\
\notag && ((\bar{a} \subseteq \{\bar{b},\bar{c},\bar{d}\}) \Leftrightarrow %
((\bar{a} \subseteq \bar{b}) \vee (\bar{a} \subseteq \bar{c}) \vee %
(\bar{a} \subseteq \bar{d}))) \wedge \\
\notag && ((\bar{a} \subseteq \{\bar{b},\bar{c},\bar{d},\bar{e}\}) \Leftrightarrow %
((\bar{a} \subseteq \bar{b}) \vee (\bar{a} \subseteq \bar{c}) \vee %
(\bar{a} \subseteq \bar{d}) \vee (\bar{a} \subseteq \bar{e}))) \wedge \\
\notag && ((\bar{a} \subseteq \{\bar{b},\bar{c},\bar{d},\bar{e},\bar{f}\}) \Leftrightarrow %
((\bar{a} \subseteq \bar{b}) \vee (\bar{a} \subseteq \bar{c}) \vee %
(\bar{a} \subseteq \bar{d}) \vee (\bar{a} \subseteq \bar{e}) \vee %
(\bar{a} \subseteq \bar{f}))) \wedge \\
\notag && ((\bar{a} \subseteq \{\bar{b},\bar{c},\bar{d},\bar{e},\bar{f},\bar{g}\}) %
\Leftrightarrow ((\bar{a} \subseteq \bar{b}) \vee (\bar{a} \subseteq \bar{c}) \vee %
(\bar{a} \subseteq \bar{d}) \vee (\bar{a} \subseteq \bar{e}) \vee (\bar{a} \subseteq \bar{f}) \\
\notag && \vee (\bar{a} \subseteq \bar{g}))) \wedge \\
\notag && ((\bar{a} \subseteq \{\bar{b},\bar{c},\bar{d},\bar{e},\bar{f},\bar{g},\bar{h}\}) %
\Leftrightarrow ((\bar{a} \subseteq \bar{b}) \vee (\bar{a} \subseteq \bar{c}) \vee %
(\bar{a} \subseteq \bar{d}) \vee (\bar{a} \subseteq \bar{e}) \vee (\bar{a} \subseteq \bar{f}) \\
\notag && \vee (\bar{a} \subseteq \bar{g}) \vee (\bar{a} \subseteq \bar{h}))) \wedge \\
\notag && ((\bar{a} \subseteq %
\{\bar{b},\bar{c},\bar{d},\bar{e},\bar{f},\bar{g},\bar{h},\bar{i}\}) %
\Leftrightarrow ((\bar{a} \subseteq \bar{b}) \vee (\bar{a} \subseteq \bar{c}) \vee %
(\bar{a} \subseteq \bar{d}) \vee (\bar{a} \subseteq \bar{e}) \vee (\bar{a} \subseteq \bar{f}) \\
\notag && \vee (\bar{a} \subseteq \bar{g}) \vee (\bar{a} \subseteq \bar{h}) %
\vee (\bar{a} \subseteq \bar{i}))) \wedge \\
\notag && ((\bar{a} \subseteq %
\{\bar{b},\bar{c},\bar{d},\bar{e},\bar{f},\bar{g},\bar{h},\bar{i},\bar{j}\}) %
\Leftrightarrow ((\bar{a} \subseteq \bar{b}) \vee (\bar{a} \subseteq \bar{c}) \vee %
(\bar{a} \subseteq \bar{d}) \vee (\bar{a} \subseteq \bar{e}) \vee \\
\notag && (\bar{a} \subseteq \bar{f}) \vee (\bar{a} \subseteq \bar{g}) \vee %
(\bar{a} \subseteq \bar{h}) \vee (\bar{a} \subseteq \bar{i}) \vee %
(\bar{a} \subseteq \bar{j}))), \\
\notag && ((\bar{a} \subseteq %
\{\bar{b},\bar{c},\bar{d},\bar{e},\bar{f},\bar{g},\bar{h},\bar{i},\bar{j},\bar{k}\}) %
\Leftrightarrow ((\bar{a} \subseteq \bar{b}) \vee (\bar{a} \subseteq \bar{c}) \vee %
(\bar{a} \subseteq \bar{d}) \vee (\bar{a} \subseteq \bar{e}) \vee \\
\notag && (\bar{a} \subseteq \bar{f}) \vee (\bar{a} \subseteq \bar{g}) \vee %
(\bar{a} \subseteq \bar{h}) \vee (\bar{a} \subseteq \bar{i}) \vee %
(\bar{a} \subseteq \bar{j}) \vee (\bar{a} \subseteq \bar{k}))), \\
\notag && ((\bar{a} \subseteq %
\{\bar{b},\bar{c},\bar{d},\bar{e},\bar{f},\bar{g},\bar{h},\bar{i},\bar{j},\bar{k},\bar{l}\}) %
\Leftrightarrow ((\bar{a} \subseteq \bar{b}) \vee (\bar{a} \subseteq \bar{c}) \vee %
(\bar{a} \subseteq \bar{d}) \vee (\bar{a} \subseteq \bar{e}) \vee \\
\notag && (\bar{a} \subseteq \bar{f}) \vee (\bar{a} \subseteq \bar{g}) \vee %
(\bar{a} \subseteq \bar{h}) \vee (\bar{a} \subseteq \bar{i}) \vee %
(\bar{a} \subseteq \bar{j}) \vee (\bar{a} \subseteq \bar{k}) \vee %
(\bar{a} \subseteq \bar{l}))), \\
\label{4.136} && (\bar{a}\bar{b}\bar{c} = \bar{d}\bar{b}\bar{e}\bar{f}\bar{g}) \wedge %
\neg(\bar{b} \subseteq \{\bar{a},\bar{c},\bar{d},\bar{e},\bar{f},\bar{g}\}) \wedge %
\neg(\bar{f} \subseteq \{\bar{a},\bar{b},\bar{c},\bar{d},\bar{e},\bar{g}\}) \wedge \\
\notag && ((\bar{b} \rightarrow \bar{h}) \| (\bar{f} \rightarrow \bar{i})) %
\Rightarrow (\bar{a}\bar{h}\bar{c} = \bar{d}\bar{h}\bar{e}\bar{i}\bar{g}), \\
\label{4.137} && (\bar{a}\bar{b}\bar{c}\bar{d}\bar{e} = \bar{f}) %
\wedge \neg(\bar{b} \subseteq \{\bar{a},\bar{c},\bar{d},\bar{e},\bar{f}\}) %
\wedge \neg(\bar{d} \subseteq \{\bar{a},\bar{b},\bar{c},\bar{e},\bar{f}\}) %
\wedge ((\bar{b} \rightarrow \bar{g}) \| (\bar{d} \rightarrow \bar{h})) \\
\notag && \Rightarrow (\bar{a}\bar{g}\bar{c}\bar{h}\bar{e} = \bar{f}), \\
\label{4.138} && a \rightarrow 1|[aba], \\
\label{4.139} && b \rightarrow +|-, \\
\label{4.140} && c|d \rightarrow e|f|g, \\
\label{4.141} && e \rightarrow +|+e, \\
\label{4.142} && f \rightarrow -|-f, \\
\label{4.143} && g \rightarrow /|/g, \\
\label{4.144} && (h \rightarrow +) \| (i \rightarrow -), \\
\label{4.145} && (h \rightarrow +h) \| (i \rightarrow -i), \\
\label{4.146} && (i \rightarrow -) \| (h \rightarrow +), \\
\label{4.147} && (i \rightarrow -i) \| (h \rightarrow +h), \\
\label{4.148} && (h \rightarrow +) \| (j \rightarrow /), \\
\label{4.149} && (h \rightarrow +h) \| (j \rightarrow /j), \\
\label{4.150} && k \rightarrow [1+1]|[1+k], \\
\label{4.151} && l \rightarrow 1|[1+l], \\
\label{4.152} && a < [1+a], \\
\label{4.153} && (\bar{a} < \bar{b}) \wedge \bar{c} \Rightarrow %
(([\bar{a}+\bar{c}] < [\bar{b}+\bar{c}]) \wedge ([\bar{c}-\bar{b}] < [\bar{c}-\bar{a}])), \\
\label{4.154} && ([1-1] < \bar{a}) \wedge (\bar{a} < \bar{b}) \wedge ([1-1] < \bar{c}) %
\Rightarrow (([\bar{a}--\bar{c}] < [\bar{b}--\bar{c}]) \wedge \\
\notag && ([\bar{c}--\bar{b}] < [\bar{c}--\bar{a}])), \\
\label{4.155} && (1 < \bar{a}) \wedge (1 < \bar{b}) \Rightarrow %
(1 < [\bar{a}--f\bar{b}]), \\
\label{4.156} && (1 < \bar{a}) \wedge (\bar{a} < \bar{b}) \Rightarrow %
(1 < [\bar{b}/g\bar{a}]), \\
\label{4.157} && (1 < \bar{a}) \wedge (\bar{a} < \bar{b}) \wedge (1 < \bar{c}) %
\Rightarrow ([\bar{a}e\bar{c}] < [\bar{b}e\bar{c}]), \\
\label{4.158} && (1 < \bar{a}) \wedge (1 < \bar{b}) \wedge (1 < \bar{c}) \wedge %
([\bar{a}e\bar{c}] < [\bar{b}e\bar{c}]) \Rightarrow (\bar{a} < \bar{b}), \\
\label{4.159} && (1 < \bar{a}) \wedge (1 < \bar{b}) \wedge (\bar{b} < \bar{c}) %
\Rightarrow ([\bar{a}e\bar{b}] < [\bar{a}e\bar{c}]), \\
\label{4.160} && (1 < \bar{a}) \wedge (1 < \bar{b}) \wedge (1 < \bar{c}) \wedge %
([\bar{a}e\bar{b}] < [\bar{a}e\bar{c}]) \Rightarrow (\bar{b} < \bar{c}), \\
\label{4.161} && \bar{a} \wedge \bar{b} \wedge \bar{c} \Rightarrow %
([\bar{a}c\bar{b}d\bar{c}] = [[\bar{a}c\bar{b}]d\bar{c}]), \\
\label{4.162} && \bar{a} \wedge \bar{b} \Rightarrow ([\bar{a}-\bar{b}] = [\bar{a}/\bar{b}]), \\
\label{4.163} && \bar{a} \wedge \neg (\bar{b} = [1-1]) \Rightarrow %
([\bar{a}--\bar{b}] = [\bar{a}//\bar{b}]), \\
\label{4.164} && \bar{a} \Rightarrow ([\bar{a}-\bar{a}] = [1-1]), \\
\label{4.165} && \bar{a} \wedge \bar{b} \Rightarrow ([\bar{a}+\bar{b}] = [\bar{b}+\bar{a}]), \\
\label{4.166} && \bar{a} \wedge \bar{b} \Rightarrow ([[\bar{a}-\bar{b}]+\bar{b}] = \bar{a}), \\
\label{4.167} && \bar{a} \wedge \bar{b} \wedge \bar{c} \Rightarrow %
([\bar{a}-\bar{b}+\bar{c}] = [\bar{a}+\bar{c}-\bar{b}]), \\
\label{4.168} && \bar{a} \wedge \bar{b} \wedge \bar{c} \Rightarrow %
([\bar{a}+[\bar{b}+\bar{c}]] = [[\bar{a}+\bar{b}]+\bar{c}]), \\
\label{4.169} && \bar{a} \wedge \bar{b} \wedge \bar{c} \Rightarrow %
([\bar{a}+[\bar{b}-\bar{c}]] = [[\bar{a}+\bar{b}]-\bar{c}]), \\
\label{4.170} && \bar{a} \wedge \bar{b} \wedge \bar{c} \Rightarrow %
([\bar{a}-[\bar{b}+\bar{c}]] = [[\bar{a}-\bar{b}]-\bar{c}]), \\
\label{4.171} && \bar{a} \wedge \bar{b} \wedge \bar{c} \Rightarrow %
([\bar{a}-[\bar{b}-\bar{c}]] = [[\bar{a}-\bar{b}]+\bar{c}]), \\
\label{4.172} && \bar{a} \Rightarrow ([\bar{a}++1] = \bar{a}), \\
\label{4.173} && \neg (\bar{a} = [1-1]) \Rightarrow ([\bar{a}--\bar{a}] = 1), \\
\label{4.174} && \bar{a} \wedge \bar{b} \Rightarrow ([\bar{a}++\bar{b}] = [\bar{b}++\bar{a}]), \\
\label{4.175} && \bar{a} \wedge \bar{b} \wedge \bar{c} \Rightarrow %
([\bar{a}++[\bar{b}++\bar{c}]] = [[\bar{a}++\bar{b}]++\bar{c}]), \\
\label{4.176} && \bar{a} \wedge \bar{b} \wedge \bar{c} \Rightarrow %
([\bar{a}++[\bar{b}+\bar{c}]] = [[\bar{a}++\bar{b}]+[\bar{a}++\bar{c}]]), \\
\label{4.177} && \bar{a} \wedge \bar{b} \wedge \bar{c} \Rightarrow %
([\bar{a}++[\bar{b}-\bar{c}]] = [[\bar{a}++\bar{b}]-[\bar{a}++\bar{c}]]), \\
\label{4.178} && \bar{a} \wedge \neg (\bar{b} = [1-1]) \Rightarrow %
([[\bar{a}--\bar{b}]++\bar{b}] = \bar{a}), \\
\label{4.179} && \bar{a} \wedge \neg (\bar{b} = [1-1]) \wedge \bar{c} %
\Rightarrow (([\bar{a}--\bar{b}++\bar{c}] = [\bar{a}++\bar{c}--\bar{b}]) \wedge \\
\notag && ([[\bar{a}+\bar{c}]--\bar{b}] = [[\bar{a}--\bar{b}]+[\bar{c}--\bar{b}]]) \wedge %
([[\bar{a}-\bar{c}]--\bar{b}] = \\
\notag && [[\bar{a}--\bar{b}]-[\bar{c}--\bar{b}]])), \\
\label{4.180} && \bar{a} \wedge \neg (\bar{b} = [1-1]) \wedge \neg (\bar{c} = [1-1]) %
\Rightarrow (([\bar{a}++[\bar{b}--\bar{c}]] = [[\bar{a}++\bar{b}]--\bar{c}]) \wedge \\
\notag && ([\bar{a}--[\bar{b}++\bar{c}]] = [[\bar{a}--\bar{b}]--\bar{c}]) \wedge %
([\bar{a}--[\bar{b}--\bar{c}]] = [[\bar{a}--\bar{b}]++\bar{c}])), \\
\label{4.181} && \bar{a} \Rightarrow ([\bar{a}+++1] = \bar{a}), \\
\label{4.182} && \bar{a} \Rightarrow ([1+++\bar{a}] = 1), \\
\label{4.183} && \neg (\bar{a} = [1-1]) \Rightarrow ([\bar{a}+++[1-1]] = 1), \\
\label{4.184} && ([1-1] < \bar{a}) \Rightarrow ([[1-1]+++\bar{a}] = [1-1]), \\
\label{4.185} && ([1-1] < \bar{a}) \wedge \neg (\bar{b} = [1-1]) \wedge \bar{c} \Rightarrow %
(([\bar{a}---\bar{b}+++\bar{b}] = \bar{a}) \wedge \\
\notag && ([\bar{a}---\bar{b}+++\bar{c}] = [\bar{a}+++\bar{c}---\bar{b}]) \wedge %
([\bar{a}+++[\bar{c}--\bar{b}]] = \\
\notag && [[\bar{a}+++\bar{c}]---\bar{b}])), \\
\label{4.186} && ([1-1] < \bar{a}) \wedge ([1-1] < \bar{b}) \wedge \bar{c} \Rightarrow %
(([\bar{a}+++[\bar{b}///\bar{a}]] = \bar{b}) \wedge \\
\notag && ([[\bar{a}+++\bar{c}]///\bar{b}] = [\bar{c}++[\bar{a}///\bar{b}]]) \wedge %
([[\bar{a}--\bar{b}]+++\bar{c}] = \\
\notag && [[\bar{a}+++\bar{c}]--[\bar{b}+++\bar{c}]])), \\
\label{4.187} && ([1-1] < \bar{a}) \wedge ([1-1] < \bar{b}) \wedge ([1-1] < \bar{c}) %
\Rightarrow (([[\bar{a}///\bar{c}]--[\bar{b}///\bar{c}]] = \\
\notag && [\bar{a}///\bar{b}]) \wedge ([[\bar{a}++\bar{b}]///\bar{c}] = %
[[\bar{a}///\bar{c}]+[\bar{b}///\bar{c}]]) \wedge ([[\bar{a}--\bar{b}]///\bar{c}] = \\
\notag && [[\bar{a}///\bar{c}]-[\bar{b}///\bar{c}]])), \\
\label{4.188} && \neg (\bar{a} = [1-1]) \wedge \neg (\bar{b} = [1-1]) \wedge %
\neg (\bar{c} = [1-1]) \Rightarrow (([[\bar{a}++\bar{b}]+++\bar{c}] = \\
\notag && [[\bar{a}+++\bar{c}]++[\bar{b}+++\bar{c}]]) \wedge %
([\bar{a}+++[\bar{b}++\bar{c}]] = [[\bar{a}+++\bar{b}]+++\bar{c}]) \\
\notag && \wedge ([\bar{a}+++[\bar{b}+\bar{c}]] = %
[[\bar{a}+++\bar{b}]++[\bar{a}+++\bar{c}]]) \wedge ([\bar{a}+++[\bar{b}-\bar{c}]] = \\
\notag && [[\bar{a}+++\bar{b}]--[\bar{a}+++\bar{c}]])), \\
\label{4.189} && ([1-1] < \bar{a}) \wedge \neg (\bar{b} = [1-1]) \wedge %
\neg (\bar{c} = [1-1]) \Rightarrow (([\bar{a}---[\bar{b}++\bar{c}]] = \\
\notag && [[\bar{a}---\bar{b}]---\bar{c}]) \wedge %
([\bar{a}---[\bar{b}--\bar{c}]] = [[\bar{a}---\bar{b}]+++\bar{c}])), \\
\label{4.190} && ([1-1] < \bar{a}) \wedge ([1-1] < \bar{b}) \wedge \neg (\bar{c} = [1-1]) %
\Rightarrow (([[\bar{a}++\bar{b}]---\bar{c}] = \\
\notag && [[\bar{a}---\bar{c}]++[\bar{b}---\bar{c}]]) \wedge %
([[\bar{a}--\bar{b}]---\bar{c}] = \\
\notag && [[\bar{a}---\bar{c}]--[\bar{b}---\bar{c}]])), \\
\label{4.191} && (1 < \bar{a}) \Rightarrow ([\bar{a}+e1] = \bar{a}), \\
\label{4.192} && (1 < \bar{a}) \Rightarrow ([\bar{a}-f1] = \bar{a}), \\
\label{4.193} && (1 < \bar{a}) \Rightarrow ([1++e\bar{a}] = 1), \\
\label{4.194} && (1 < \bar{a}) \Rightarrow ([1--f\bar{a}] = 1), \\
\label{4.195} && (1 < \bar{a}) \Rightarrow ([1///\bar{a}] = [1-1]), \\
\label{4.196} && (1 < \bar{a}) \Rightarrow ([\bar{a}/g\bar{a}] = 1), \\
\label{4.197} && (1 < \bar{a}) \wedge (1 < \bar{b}) \Rightarrow %
([[\bar{a}i\bar{b}]h\bar{b}] = \bar{a}), \\
\label{4.198} && (1 < \bar{a}) \wedge (1 < \bar{b}) \Rightarrow %
([[\bar{a}h\bar{b}]i\bar{b}] = \bar{a}), \\
\label{4.199} && (1 < \bar{a}) \wedge (1 < \bar{b}) \Rightarrow %
([\bar{b}h[\bar{a}j\bar{b}]] = \bar{a}), \\
\label{4.200} && (1 < \bar{a}) \wedge (1 < \bar{b}) \Rightarrow %
([[\bar{b}h\bar{a}]j\bar{b}] = \bar{a}), \\
\label{4.201} && (1 < \bar{a}) \wedge (1 < \bar{b}) \Rightarrow %
([\bar{a}+e\bar{b}] = [\bar{a}e[\bar{a}+e[\bar{b}-1]]]), \\
\label{4.202} && \top\underline{\ \ }1\underline{\ \ }1\underline{\ \ } = [1-1], \\
\label{4.203} && \top\underline{\ \ }1\underline{\ \ }[1+1]\underline{\ \ } = 1, \\
\label{4.204} && \bot\underline{\ \ }1\underline{\ \ }1\underline{\ \ } = 1, \\
\label{4.205} && \bot\underline{\ \ }1\underline{\ \ }[1+1]\underline{\ \ } = 1, \\
\label{4.206} && \top\underline{\ \ }k\underline{\ \ }[[[1+1]++l]-1]\underline{\ \ } = %
\top\underline{\ \ }[k-1]\underline{\ \ }l\underline{\ \ }, \\
\label{4.207} && \bot\underline{\ \ }k\underline{\ \ }[[[1+1]++l]-1]\underline{\ \ } = %
\bot\underline{\ \ }[k-1]\underline{\ \ }l\underline{\ \ }, \\
\label{4.208} && \top\underline{\ \ }k\underline{\ \ }[[1+1]++l]\underline{\ \ } = %
[\top\underline{\ \ }[k-1]\underline{\ \ }l\underline{\ \ }+ %
\top\underline{\ \ }[k-1]\underline{\ \ }[l+1]\underline{\ \ }], \\
\label{4.209} && \bot\underline{\ \ }k\underline{\ \ }[[1+1]++l]\underline{\ \ } = %
[\bot\underline{\ \ }[k-1]\underline{\ \ }l\underline{\ \ }+ %
\bot\underline{\ \ }[k-1]\underline{\ \ }[l+1]\underline{\ \ }], \\
\label{4.210} && (1 < \bar{a}) \wedge %
(\top\underline{\ \ }\bar{b}\underline{\ \ }\bar{c}\underline{\ \ }) \wedge  %
(\bot\underline{\ \ }\bar{b}\underline{\ \ }\bar{c}\underline{\ \ }) \Rightarrow %
([\bar{a}+h[\top\underline{\ \ }\bar{b}\underline{\ \ }\bar{c}\underline{\ \ }-- %
\bot\underline{\ \ }\bar{b}\underline{\ \ }\bar{c}\underline{\ \ }]] = \\
\notag && [[\bar{a}+h\top\underline{\ \ }\bar{b}\underline{\ \ }\bar{c}\underline{\ \ }]-i %
\bot\underline{\ \ }\bar{b}\underline{\ \ }\bar{c}\underline{\ \ }]) \\
\notag && \}. %
\end{eqnarray}
\end{definition}

It should be noted that (\ref{4.202}) $ \sim $ (\ref{4.210}) restrict %
$ [\top\underline{\ \ }\bar{b}\underline{\ \ }\bar{c}\underline{\ \ }-- %
\bot\underline{\ \ }\bar{b}\underline{\ \ }\bar{c}\underline{\ \ }] $ to be a Farey fraction. %
In the following, we will deduce some numbers and equalities as examples. %
\begin{eqnarray*}
(A1) & (a \rightarrow 1|[aba]) \Rightarrow (a \rightarrow 1) & by (\ref{4.138}),(\ref{4.18}) \\
(A2) & \Rightarrow (a \rightarrow [aba]) & by (\ref{4.138}),(\ref{4.18}) \\
(A3) & \Rightarrow (a \rightarrow [[aba]ba]) & by (A2),(\ref{4.17}) \\
(A4) & \Rightarrow (a \rightarrow [[1-1]-1]) & by (A1),(\ref{4.139}) \\
(A5) & (a < [1+a]) \Rightarrow (1 < [1+1]) & by (\ref{4.152}),(A1),(\ref{4.29}) \\
(A6) & \Rightarrow 1 & by (\ref{4.23}) \\
(A7) & \Rightarrow [1+1] & \ by (A5),(\ref{4.23}) \\
(A8) & \Rightarrow ([aba] < [1+[aba]]) & by (A5),(A2),(\ref{4.29}) \\
(A9) & \Rightarrow ([1+1] < [1+[1+1]]) & by (A1),(\ref{4.139}),(\ref{4.29}) \\
(A10) & \Rightarrow ([1+1] < [[1+1]+1]) & by (\ref{4.165}),(\ref{4.24}) \\
(A11) & \Rightarrow ([1+1] < [1+1+1]) & by (\ref{4.161}) \\
(A12) & \Rightarrow ([1-1] < [1+[1-1]]) & by (A8),(A1),(\ref{4.139}),(\ref{4.29}) \\
(A13) & \Rightarrow ([1-1] < [[1-1]+1]) & by (\ref{4.165}),(\ref{4.24}) \\
(A14) & \Rightarrow ([1-1] < [1-1+1]) & by (\ref{4.161}) \\
(A15) & \Rightarrow ([1-1] < 1) & by (\ref{4.166}) \\
(A16) & \Rightarrow ([1-1] < [1+1+1]) & by (A15),(A5),(A11),(\ref{4.22}) \\
(A17) & \Rightarrow [1-1] & by (A15),(\ref{4.23}) \\
(A18) & \Rightarrow [1+1+1] & by (A16),(\ref{4.23}) \\
(A19) & \Rightarrow ([1-1] < [1+1]) & by (A15),(A5),(\ref{4.22}) \\
(A20) & ([[1+1]--[1+1]] < [[1+1+1]--[1+1]]) & by (A5),(A11),(\ref{4.154}) \\
(A21) & \Rightarrow (1 < [[1+1+1]--[1+1]]) & by (A20),(\ref{4.173}) \\
(A22) & \Rightarrow [[1+1+1]--[1+1]] & by (A21),(\ref{4.23}) \\
\vdots & \vdots & \vdots
\end{eqnarray*}

Then we deduce the numbers from $ R \{ \Phi, \Psi \} $: %
\begin{eqnarray*}
& [1-1], [[1-1]-[[1+1]----[1+1]]], [[1+1+1]----[1+1]] \cdots & %
\end{eqnarray*}
\begin{eqnarray*}
(B1) & (\top\underline{\ \ }k\underline{\ \ }[[1+1]++l]\underline{\ \ } = %
[\top\underline{\ \ }[k-1]\underline{\ \ }l\underline{\ \ }+ %
\top\underline{\ \ }[k-1]\underline{\ \ }[l+1]\underline{\ \ }]) & by (\ref{4.208}) \\
(B2) & \Rightarrow (\top\underline{\ \ }[1+1]\underline{\ \ }[[1+1]++l]\underline{\ \ } = %
[\top\underline{\ \ }[[1+1]-1]\underline{\ \ }l\underline{\ \ }+ & \\
& \top\underline{\ \ }[[1+1]-1]\underline{\ \ }[l+1]\underline{\ \ }]) & by (\ref{4.150}), %
(\ref{4.18}), \\
&& (\ref{4.42}) \\
(B3) & \Rightarrow (\top\underline{\ \ }[1+1]\underline{\ \ }[[1+1]++1]\underline{\ \ } = %
[\top\underline{\ \ }[[1+1]-1]\underline{\ \ }1\underline{\ \ }+ & \\
& \top\underline{\ \ }[[1+1]-1]\underline{\ \ }[1+1]\underline{\ \ }]) & by (\ref{4.151}), %
(\ref{4.18}), \\
&& (\ref{4.42}) \\
(B4) & \Rightarrow (\top\underline{\ \ }[1+1]\underline{\ \ }[[1+1]++1]\underline{\ \ } = %
[\top\underline{\ \ }[1+1-1]\underline{\ \ }1\underline{\ \ }+ & \\
& \top\underline{\ \ }[1+1-1]\underline{\ \ }[1+1]\underline{\ \ }]) & by (\ref{4.161}) \\
(B5) & \Rightarrow (\top\underline{\ \ }[1+1]\underline{\ \ }[[1+1]++1]\underline{\ \ } = %
[\top\underline{\ \ }[1-1+1]\underline{\ \ }1\underline{\ \ }+ & \\
& \top\underline{\ \ }[1-1+1]\underline{\ \ }[1+1]\underline{\ \ }]) & by (\ref{4.167}) \\
(B6) & \Rightarrow (\top\underline{\ \ }[1+1]\underline{\ \ }[[1+1]++1]\underline{\ \ } = %
[\top\underline{\ \ }[[1-1]+1]\underline{\ \ }1\underline{\ \ }+ & \\
& \top\underline{\ \ }[[1-1]+1]\underline{\ \ }[1+1]\underline{\ \ }]) & by (\ref{4.161}) \\
(B7) & \Rightarrow (\top\underline{\ \ }[1+1]\underline{\ \ }[[1+1]++1]\underline{\ \ } = %
[\top\underline{\ \ }1\underline{\ \ }1\underline{\ \ }+ %
\top\underline{\ \ }1\underline{\ \ }[1+1]\underline{\ \ }]) & by (\ref{4.166}) \\
(B8) & \Rightarrow (\top\underline{\ \ }[1+1]\underline{\ \ }[1+1]\underline{\ \ } = %
[\top\underline{\ \ }1\underline{\ \ }1\underline{\ \ }+ %
\top\underline{\ \ }1\underline{\ \ }[1+1]\underline{\ \ }]) & by (\ref{4.172}) \\
(B9) & \Rightarrow (\top\underline{\ \ }[1+1]\underline{\ \ }[1+1]\underline{\ \ } = %
[[1-1]+1]) & by (\ref{4.202}),(\ref{4.203}) \\
(B10) & \Rightarrow (\top\underline{\ \ }[1+1]\underline{\ \ }[1+1]\underline{\ \ } = 1) %
& by (\ref{4.166}) \\
(B11) & \Rightarrow ([1-1] < \top\underline{\ \ }[1+1]\underline{\ \ }[1+1]\underline{\ \ }) %
& by (A15),(\ref{4.24}) \\
(B12) & \Rightarrow (\top\underline{\ \ }[1+1]\underline{\ \ }[1+1]\underline{\ \ }) %
& by (\ref{4.23}) \\
(B13) & (\bot\underline{\ \ }k\underline{\ \ }[[1+1]++l]\underline{\ \ } = %
[\bot\underline{\ \ }[k-1]\underline{\ \ }l\underline{\ \ }+ %
\bot\underline{\ \ }[k-1]\underline{\ \ }[l+1]\underline{\ \ }]) & by (\ref{4.209}) \\
(B14) & \Rightarrow (\bot\underline{\ \ }[1+1]\underline{\ \ }[[1+1]++l]\underline{\ \ } = %
[\bot\underline{\ \ }[[1+1]-1]\underline{\ \ }l\underline{\ \ }+ & \\
& \bot\underline{\ \ }[[1+1]-1]\underline{\ \ }[l+1]\underline{\ \ }]) & by (\ref{4.150}), %
(\ref{4.18}), \\
&& (\ref{4.42}) \\
(B15) & \Rightarrow (\bot\underline{\ \ }[1+1]\underline{\ \ }[[1+1]++1]\underline{\ \ } = %
[\bot\underline{\ \ }[[1+1]-1]\underline{\ \ }1\underline{\ \ }+ & \\
& \bot\underline{\ \ }[[1+1]-1]\underline{\ \ }[1+1]\underline{\ \ }]) & by (\ref{4.151}), %
(\ref{4.18}), \\
&& (\ref{4.42}) \\
(B16) & \Rightarrow (\bot\underline{\ \ }[1+1]\underline{\ \ }[[1+1]++1]\underline{\ \ } = %
[\bot\underline{\ \ }[1+1-1]\underline{\ \ }1\underline{\ \ }+ & \\
& \bot\underline{\ \ }[1+1-1]\underline{\ \ }[1+1]\underline{\ \ }]) & by (\ref{4.161}) \\
(B17) & \Rightarrow (\bot\underline{\ \ }[1+1]\underline{\ \ }[[1+1]++1]\underline{\ \ } = %
[\bot\underline{\ \ }[1-1+1]\underline{\ \ }1\underline{\ \ }+ & \\
& \bot\underline{\ \ }[1-1+1]\underline{\ \ }[1+1]\underline{\ \ }]) & by (\ref{4.167}) \\
(B18) & \Rightarrow (\bot\underline{\ \ }[1+1]\underline{\ \ }[[1+1]++1]\underline{\ \ } = %
[\bot\underline{\ \ }[[1-1]+1]\underline{\ \ }1\underline{\ \ }+ & \\
& \bot\underline{\ \ }[[1-1]+1]\underline{\ \ }[1+1]\underline{\ \ }]) & by (\ref{4.161}) \\
(B19) & \Rightarrow (\bot\underline{\ \ }[1+1]\underline{\ \ }[[1+1]++1]\underline{\ \ } = %
[\bot\underline{\ \ }1\underline{\ \ }1\underline{\ \ }+ %
\bot\underline{\ \ }1\underline{\ \ }[1+1]\underline{\ \ }]) & by (\ref{4.166}) \\
(B20) & \Rightarrow (\bot\underline{\ \ }[1+1]\underline{\ \ }[1+1]\underline{\ \ } = %
[\bot\underline{\ \ }1\underline{\ \ }1\underline{\ \ }+ %
\bot\underline{\ \ }1\underline{\ \ }[1+1]\underline{\ \ }]) & by (\ref{4.172}) \\
(B21) & \Rightarrow (\bot\underline{\ \ }[1+1]\underline{\ \ }[1+1]\underline{\ \ } = [1+1]) %
& by (\ref{4.204}),(\ref{4.205}) \\
(B22) & \Rightarrow ([1-1] < \bot\underline{\ \ }[1+1]\underline{\ \ }[1+1]\underline{\ \ }) %
& by (A19),(\ref{4.24}) \\
(B23) & \Rightarrow (\bot\underline{\ \ }[1+1]\underline{\ \ }[1+1]\underline{\ \ }) %
& by (\ref{4.23}) \\
(B24) & ([[1+1]+e[[1+1+1]--[1+1]]] = & \\
& [[1+1]e[[1+1]+e[[[1+1+1]--[1+1]]-1]]]) & by (A5),(A21), \\
&& (\ref{4.201}) \\
(B25) & \Rightarrow ([[1+1]++e[[1+1+1]--[1+1]]] = & \\
& [[1+1]+e[[1+1]++e[[[1+1+1]--[1+1]]-1]]]) & by (\ref{4.141}) \\
(B26) & \Rightarrow ([[1+1]+++e[[1+1+1]--[1+1]]] = & \\
& [[1+1]++e[[1+1]+++e[[[1+1+1]--[1+1]]-1]]]) & by (\ref{4.141}) \\
(B27) & \Rightarrow ([[1+1]++++[[1+1+1]--[1+1]]] = & \\
& [[1+1]+++[[1+1]++++[[[1+1+1]--[1+1]]-1]]]) & by (\ref{4.141}) \\
(B28) & \Rightarrow ([[1+1]++++[[1+1+1]--[1+1]]] = & \\
& [[1+1]+++[[1+1]++++[[[[1+1]+1]--[1+1]]-1]]]) & by (\ref{4.161}) \\
(B29) & \Rightarrow ([[1+1]++++[[1+1+1]--[1+1]]] = & \\
& [[1+1]+++[[1+1]++++ & \\
& [[[[1+1]--[1+1]]+[1--[1+1]]]-1]]]) & by (\ref{4.179}) \\
(B30) & \Rightarrow ([[1+1]++++[[1+1+1]--[1+1]]] = & \\
& [[1+1]+++[[1+1]++++[[1+[1--[1+1]]]-1]]]) & by (A19),(\ref{4.173}) \\
(B31) & \Rightarrow ([[1+1]++++[[1+1+1]--[1+1]]] = & \\
& [[1+1]+++[[1+1]++++[[[1--[1+1]]+1]-1]]]) & by (\ref{4.165}) \\
(B32) & \Rightarrow ([[1+1]++++[[1+1+1]--[1+1]]] = & \\
& [[1+1]+++[[1+1]++++[[1--[1+1]]+1-1]]]) & by (\ref{4.161}) \\
(B33) & \Rightarrow ([[1+1]++++[[1+1+1]--[1+1]]] = & \\
& [[1+1]+++[[1+1]++++[[1--[1+1]]-1+1]]]) & by (\ref{4.167}) \\
(B34) & \Rightarrow ([[1+1]++++[[1+1+1]--[1+1]]] = & \\
& [[1+1]+++[[1+1]++++[[[1--[1+1]]-1]+1]]]) & by (\ref{4.161}) \\
(B35) & \Rightarrow ([[1+1]++++[[1+1+1]--[1+1]]] = & \\
& [[1+1]+++[[1+1]++++[1--[1+1]]]]) & by (\ref{4.166}) \\
(B36) & ([[1+1]+h[\top\underline{\ \ }[1+1]\underline{\ \ }[1+1]\underline{\ \ }-- %
\bot\underline{\ \ }[1+1]\underline{\ \ }[1+1]\underline{\ \ }]] = & \\
& [[[1+1]+h\top\underline{\ \ }[1+1]\underline{\ \ }[1+1]\underline{\ \ }]-i %
\bot\underline{\ \ }[1+1]\underline{\ \ }[1+1]\underline{\ \ }]) & by (A5),(B12), \\
&& (B23),(\ref{4.210}) \\
(B37) & \Rightarrow ([[1+1]+h[1--[1+1]]] = [[[1+1]+h1]-i[1+1]]) & by (B10),(B21), \\
&& (\ref{4.41}) \\
(B38) & \Rightarrow ([[1+1]++h[1--[1+1]]] = & \\
& [[[1+1]++h1]--i[1+1]]) & by (\ref{4.136}),(\ref{4.145}) \\
(B39) & \Rightarrow ([[1+1]+++h[1--[1+1]]] = & \\
& [[[1+1]+++h1]---i[1+1]]) & by (\ref{4.136}),(\ref{4.145}) \\
(B40) & \Rightarrow ([[1+1]++++[1--[1+1]]] = & \\
& [[[1+1]++++1]----[1+1]]) & by (\ref{4.136}),(\ref{4.144}) \\
(B41) & ([[1+1]+e1] = [1+1]) & by (A5),(\ref{4.191}) \\
(B42) & \Rightarrow ([[1+1]++e1] = [1+1]) & by (\ref{4.141}) \\
(B43) & \Rightarrow ([[1+1]+++e1] = [1+1]) & by (\ref{4.141}) \\
(B44) & \Rightarrow ([[1+1]++++1] = [1+1]) & by (\ref{4.141}) \\
(B45) & \Rightarrow ([[1+1]++++[1--[1+1]]] = [[1+1]----[1+1]]) & by (B40),(B44), \\
&& (\ref{4.37}) \\
(B46) & \Rightarrow ([[1+1]++++[[1+1+1]--[1+1]]] = & \\
& [[1+1]+++[[1+1]----[1+1]]]) & by (B35),(B45), \\
&& (\ref{4.37}) \\
\vdots & \vdots & \vdots
\end{eqnarray*}

Then we deduce the equalities on deducible numbers from $ R \{ \Phi, \Psi \} $: %
\begin{eqnarray*}
& [[1+1]++[[1+1]---[1+1]]] = [[[1+1]---[1+1]]++[1+1]], & \\
& [[1+1]++++[[1+1+1]--[1+1]]] = [[1+1]+++[[1+1]----[1+1]]] & \\
& \vdots & %
\end{eqnarray*}

The deducible numbers correspond to the real numbers as follows: %
\begin{eqnarray*}
\vdots & \vdots & \vdots, \\
\ [[1-1]-[[1+1+1]----[1+1]]] & \equiv & \\
\vdots & \vdots & \vdots, \\
\ [[1-1]-[[1+1]---[1+1]]] & \equiv & -\sqrt[2]{2}, \\
\vdots & \vdots & \vdots, \\
\ [[1-1]-[[1+1]---[1+1+1]]] & \equiv & -\sqrt[3]{2}, \\
\vdots & \vdots & \vdots, \\
\ [[1-1]-1] & \equiv & -1, \\
\vdots & \vdots & \vdots, \\
\ [[1-1]-[1--[1+1]]] & \equiv & -\frac{1}{2}, \\
\vdots & \vdots & \vdots, \\
\ [1-1] & \equiv & 0, \\
\vdots & \vdots & \vdots, \\
\ [1-[1--[1+1]]] & \equiv & \frac{1}{2}, \\
\vdots & \vdots & \vdots, \\
\ [[1+1]///[1+1+1]] & \equiv & \log_{3} 2, \\
\vdots & \vdots & \vdots, \\
\ 1 & \equiv & 1, \\
\vdots & \vdots & \vdots, \\
\ [[1+1]---[1+1+1]] & \equiv & \sqrt[3]{2}, \\
\vdots & \vdots & \vdots, \\
\ [[1+1]---[1+1]] & \equiv & \sqrt[2]{2}, \\
\vdots & \vdots & \vdots, \\
\ [1+[1--[1+1]]] & \equiv & \frac{3}{2}, \\
\vdots & \vdots & \vdots, \\
\ [[1+1+1]///[1+1]] & \equiv & \log_{2} 3, \\
\vdots & \vdots & \vdots, \\
\ [[1+1+1]----[1+1]] & \equiv & \\
\vdots & \vdots & \vdots, \\
\ [1+1] & \equiv & 2, \\
\vdots & \vdots & \vdots, \\
\ [1+1+[1--[1+1]]] & \equiv & \frac{5}{2}, \\
\vdots & \vdots & \vdots, \\
\ [1+1+1] & \equiv & 3, \\
\vdots & \vdots & \vdots, \\
\ [[[1+1+1]----[1+1]]++[1+1]] & \equiv & \\
\vdots & \vdots & \vdots .
\end{eqnarray*}

The equalities on deducible numbers correspond to addition, subtraction, multiplication, %
division, power operation and more other operations in real number system. Note that in %
the correspondence above some irrational numbers such as $ [[1+1+1]----[1+1]] $, %
$ [[[1+1+1]----[1+1]]++[1+1]] $ and $ [[1-1]-[[1+1+1]----[1+1]]] $ do not correspond to any %
irrational number based on traditional operations, which however can be constructed by the %
logical calculus $ R \{ \Phi, \Psi \} $. %

Then the logical calculus $ R \{ \Phi, \Psi \} $ not only derives the irrational numbers, but %
also makes its deducible numbers join in algebraical operations. So the logical calculus %
$ R \{ \Phi, \Psi \} $ intuitively and logically denote real number system. %

%%%%%%%%%%%%%%%%%%%%%%%%%%%%%%%%%%%%%%%%%%%%%%%%%%%%%%%%%%%%%

% BibTeX users please use
% \bibliographystyle{}
% \bibliography{}

\begin{thebibliography}{}
\bibitem{Ref1}
    G. H. Hardy and E. M. Wright, An Introduction to the Theory of
    Numbers, Fifth Edition, Oxford University Press, 1979.

\bibitem{Ref2}
    Edward A. Bender and S. Gill Williamson, A Short Course in Discrete
    Mathematics, Dover Publications, 2005.

\bibitem{Ref3}
    Walter Rudin, Principles of Mathematical Analysis, Third
    Edition, McGraw-Hill, 1976.

\bibitem{Ref4}
    Victor J. Katz, A History of Mathematics: Brief Edition,
    Pearson Education, Inc., 2004.

\bibitem{Ref5}
    Paul R. Halmos, Naive Set Theory, Springer-Verlag New York, Inc., 1974.

\bibitem{Ref6}
    H. Garth Dales and W. Hugh Woodin, Super-Real Fields, Clarendon Press, 1996.

\bibitem{Ref7}
    John H. Conway, On Numbers and Games, Academic Press, 1976.

\bibitem{Ref8}
    Herbert B. Enderton, Elements of Set Theory, Academic Press,
    New York, 1977.

\bibitem{Ref9}
    Peter Linz, An Introduction to Formal Languages and Automata, Third
    Edition, Jones and Bartlett Publishers, Inc., 2001.

\bibitem{Ref10}
    Anil Nerode and Richare A. Shore, Logic for Application, Second
    Edition, Springer-Verlag New York, Inc., 1997.

\bibitem{Ref11}
    Herbert B. Enderton, A Mathematical Introduction to Logic, Second
    Edition, Elsevier, 2001.

\bibitem{Ref12}
    H.-D. Ebbinghaus, J. Flum, W. Thomas, Mathematical Logic, Second
    Edition, Springer Science+Business Media, Inc., 1994.

\bibitem{Ref13}
    C. Haros, Tables pour evaluer une fraction ordinaire avec autand de
    decimals qu'on voudra; et pour trover la fraction ordinaire la plus
    simple, et qui a approche sensiblement d'une fraction decimale, J.
    Ecole Polytechn. 4 (1802), 364-368.
\end{thebibliography}
%
% Non-BibTeX users please use

\end{document}